\documentclass[11pt]{amsart}%
\usepackage{amsmath}
\usepackage{amssymb}
\usepackage{amsthm}
\usepackage[mathscr]{eucal}
\usepackage{mathrsfs}
\usepackage{psfrag}
\usepackage{fullpage}
\usepackage{color}
\usepackage{eucal}
\usepackage[dvips]{graphicx}
\usepackage{amsfonts}%
\usepackage{amsaddr}
\providecommand{\U}[1]{\protect\rule{.1in}{.1in}}
\newtheorem{proposition}{Proposition}[section]
\newtheorem{theorem}[proposition]{Theorem}

\newtheorem{lemma}[proposition]{Lemma}

\newtheorem{definition}[proposition]{Definition}
\newtheorem{remark}[proposition]{Remark}
\newtheorem{example}[proposition]{Example}

\numberwithin{equation}{section}
\numberwithin{proposition}{section}

\numberwithin{equation}{section}
\numberwithin{proposition}{section}

\begin{document}
\title{Irreversible Langevin samplers and variance reduction: a large deviations approach}

\author{Luc Rey-Bellet}
\address{Department of  Mathematics and Statistics\\
University of Massachusetts Amherst, Amherst, MA, 01003}
\email{luc@math.umass.edu}

\author{Konstantinos Spiliopoulos}
\address{Department of  Mathematics and Statistics\\
Boston University, Boston, MA, 02215}
\email{kspiliop@math.bu.edu}
\thanks{K.S. was partially supported by the National Science Foundation
(DMS 1312124). LRB was partially supported by the National Science Foundation (DMS 1109316) and the Department of Energy-ASCR
(ER 26161) and he also thanks Natesh Pillai for bringing the problem to his attention.}

\date{\today}

\begin{abstract}
In order to sample from a given target distribution (often of Gibbs type), the Monte Carlo Markov chain method consists
in constructing an ergodic Markov process whose invariant measure is the target distribution.
By sampling the Markov process one can then compute,  approximately, expectations of observables with respect to the target
distribution.  Often the Markov processes used in practice are time-reversible (i.e., they satisfy detailed balance), but our main
goal here is to assess and quantify how the addition of a non-reversible part to the process can be used to improve the
sampling properties.  We focus on the diffusion setting (overdamped Langevin equations) where the drift consists of a gradient
vector field as well as another drift which breaks the reversibility of the process but  is chosen to preserve
the Gibbs measure.  In this paper we use the large deviation rate function for the empirical measure as a tool to analyze
the speed of convergence to the invariant measure.
We show that the addition of an irreversible drift leads to a larger rate function and it strictly improves the speed of convergence of
ergodic average for (generic smooth) observables. We also deduce from this result that the asymptotic variance
decreases under the addition of the irreversible drift and we give an explicit characterization of the observables whose
variance is not reduced reduced, in terms of a nonlinear Poisson equation. Our theoretical results are illustrated and supplemented by numerical simulations.
\end{abstract}

\maketitle

\textbf{AMS}: 60F05, 60F10, 60J25, 60J60, 65C05, 82B80

\section{Introduction}\label{S:Introduction}

In a wide range of applications it is often of interest to sample from a given high-dimensional distribution.  However, often,
the target distribution, say $\bar{\pi}(dx)$, is known only up to normalizing constants and then one has to rely on approximations.
In practice, one often relies on approximations using Markov processes that have the particular target distributions as their
invariant measure, as for example in Monte Carlo Markov Chain methods. Closely related, in steady-state simulations one is often
interested in quantities of the form $\int_{E}f(x)\bar{\pi}(dx)$, where $E$ is the state space and $f$ is a given function.
When closed-form  evaluation of such integrals is prohibitive, one considers a Markov process $X_{t}$ which has $\bar{\pi}$ as
its invariant distribution  and  under the assumption that $X_{t}$ is positive recurrent, the ergodic theorem gives
\begin{equation}
\frac{1}{t}\int_{0}^{t}f(X_{s})ds\rightarrow \int_{E}f(x)\bar{\pi}(dx), \text{ a.s. as }t\rightarrow\infty \,,\label{Eq:ErgodicTheorem}
\end{equation}
for all $f \in L^1(\bar{\pi})$.
Hence, the estimator $f_{t}\equiv\frac{1}{t}\int_{0}^{t}f(X_{s})ds$  can be used to approximate the expectation
${\bar f} \equiv \int_{E}f(x)\bar{\pi}(dx)$.

Standard criteria to analyze the degree of efficiency of a simulation method relies on the ergodic properties of the Markov process.
The spectral gap of the semigroup in $L^2(\pi)$ (or in other functional settings), which provides a bound for the distance
between the distribution of $X_t$ and $\pi$, as well as the asymptotic variance of $f_t$ are commonly used,
see for example \cite{AmitGrenander1991,AtheyaDossSethuraman1996,BedardRosenthal2008,Bierkens,DiaconisHolmesNeal2010,DiaconisMiclo2013,FrankeHwangPaiSheu2010,FrigessiHwangSheuStefano1993,FrigessiHwangYounes1992,GilksRoberts1996,HwangMaSheu1993,HwangMaSheu2005,LelievreNierPavliotis2012,Mira2001b,MiraGeyer2000,Neal2004,MergessenTweedie1996,Peskun,RobertsRosenthal2004,SunGomezSchmidhuber,Tierney}.
A couple of years ago, in \cite{DupuisDoll1,DupuisDoll2},  the theory of large deviations, specifically the rate function for the
empirical measure, has been proposed as a comparison tool to assess Monte-Carlo methods and used to analyze
the swapping algorithm. In this paper we use this criterium as a guide to design and analyze non-reversible Markov processes and
compare them  with reversible ones. We show that the rate function increases under the
addition of an irreversible drift. This is shown to improve the convergence properties of
the ergodic average $f_t$ for generic (smooth) observables.  We prove as well that a fine analysis of the
large deviation rate function  allows us to show that the asymptotic variance for generic smooth observables decreases.

In this paper, we specialize to the diffusion setting: to sample the Gibbs measure ${\bar{\pi}}$
on the set $E$ with density
\begin{equation*}
\frac{e^{-2 U(x)}}{\int_E e^{- 2 U(x)}dx}
\end{equation*}
one can consider the (time-reversible) Langevin equation
\begin{equation}\label{reversibleS}
dX_{t}=-\nabla U(X_{t})dt+  dW_{t}\,,
\end{equation}
whose invariant measure is ${\bar{\pi}}$.
There are however many other stochastic differential equations with the same invariant measure, for example
the family of equations
\begin{equation}\label{irreversibleS}
dX_{t}=\left[-\nabla U(X_{t})+C(X_{t})\right]dt+dW_{t}\,,
\end{equation}
where the vector field $C(x)$ satisfies the condition
\[
\textrm{div}(C e^{-2U}) \,=\,0.
\]

This constraint ensures that ${\bar \pi}$ remains the unique invariant measure, but then the Markov process is time-reversible only if $C=0$. There are many possible choices for the vector field $C(x)$. Indeed,
since $\textrm{div}(C e^{-2U})=0$ is equivalent to
\[
\textrm {div}(C) = 2 C \nabla U\,,
\]
we can choose for example $C$ to be both divergence free and orthogonal to $\nabla U$. In any dimension one can for example
set $C=S\nabla U$ where $S$ is an (arbitrary) anti-symmetric matrix $S$. More generally,
by Theorem 5.3 of \cite{Barbarosie2011}  any divergence free vector field in dimension $d$ can be written, locally,
as the exterior (or wedge) product $ C = \nabla V_1 \wedge \cdots \wedge \nabla V_{n-1}$ for some $V_i \in \mathcal{C}^1(E)$.  Therefore for our purpose we can pick $C$ of the form
\[
C = \nabla U \wedge \nabla V_2 \cdots \wedge \nabla V_{n-1}\,.
\]
for arbitrary $V_2, \cdots, V_{n-1}\in \mathcal{C}^1(E)$, and this guarantees that $C\nabla U=0$ by the properties of the exterior
product.

The main result in \cite{HwangMaSheu2005} is that the absolute value of the second largest eigenvalue  of the Markov
semigroup in $L^2(\bar{\pi})$  strictly decreases under a natural non-degeneracy condition on $C$
(the corresponding eigenspace should not be invariant under the action of the added drift $C$).
More detailed results on the  spectral gap are in \cite{ConstantinKiselevRyshikZlatos2008,FrankeHwangPaiSheu2010} where
the authors consider diffusions on compact manifolds  with $U=0$ and a one-parameter families of perturbations $C = \delta C_0$ for
$\delta \in \mathbb{R}$ and $C_0$ is some divergence vector field. In these papers the behavior of the spectral gap is related to the
ergodic properties of the flow generated by $C$
(for example if the flow is weak-mixing then the second largest  eigenvalue tends to $0$ as $\delta \to \infty$). Further, a
detailed analysis of linear diffusion processes with $U(x)= \frac{1}{2}x^{T}Ax$ and $C= J Ax$ for a antisymmetric $J$ can be found
in \cite{HwangMaSheu1993,LelievreNierPavliotis2012}  where the optimal  choice of $J$ is  determined.

We consider here the same class of problems but we take the large deviations rate function as a
measure of  the speed of convergence to equilibrium and deduce from it results on the asymptotic variance for a given observable.
While the spectral gap measures the distance of the distribution of $X_t$ compared to the invariant distribution, from a practical
Monte-Carlo point of view one is often more interested in the distribution of the ergodic average $t^{-1} \int_0^t f(X_s) \, ds$ and
how likely it  is that this average differs from the average $\int f d \bar{\pi}$.
It will be useful to consider in a  first step the empirical measure
\begin{equation}\label{emp-measure}
\pi_{t} \equiv \frac{1}{t}\int_{0}^{t} \delta _{X_s} \, ds
\end{equation}
which converges to $\bar{\pi}$ almost surely.  Let us assume that we have a large deviation principle for the family of measures $\pi_t$, which we write, symbolically as
\begin{equation*}
\mathbb{P} \left\{ \pi_{t}  \approx \mu \right\} \asymp e^{- t I_C(\mu)}.
\end{equation*}

Here $\asymp$ denotes logarithmic equivalence (the formal definition is given in Definition \ref{Def:LDP}). Then, the rate function
$I_C(\mu)$ which is non-negative and vanishes if and only if $\mu=\bar{\pi}$ quantifies the exponential rate at which the random
measure $\pi_t$ converges to $\bar{\pi}$.  Clearly, the larger $I_C$ is, the faster
the convergence occurs.

Breaking detailed balance has been shown to accelerate convergence to equilibrium for Markov chains by increasing spectral gap and/or
decreasing asymptotic variance and for diffusions by increasing spectral gap, e.g.,
\cite{Bierkens,DiaconisHolmesNeal2010,DiaconisMiclo2013,FrankeHwangPaiSheu2010,FrigessiHwangSheuStefano1993,FrigessiHwangYounes1992,HwangMaSheu1993,HwangMaSheu2005,Mira2001b,MiraGeyer2000,Neal2004,SunGomezSchmidhuber}.
The novelty of the present paper lies in that (a): we use large deviations theory in a novel way to characterize convergence to equilibrium,
(b): we prove that asymptotic variance is also decreased when breaking detailed balance for diffusions, and (c): we derive a Poisson
equation which characterizes when irreversible perturbations lead to strict improvement in performance.

Our first key result here is that if $\mu(dx) =p(x) dx$ has a smooth density $p$ and satisfies the non-degeneracy condition
${\rm div}(pC)  \not=0$, the large deviation rate function strictly increases, $I_C(\mu)> I_0(\mu)$, when one adds a non-zero
appropriate  drift  $C(x)$ to make the process $X_{t}$ irreversible,  see Theorem \ref{T:measure1}.
Moreover, specializing to perturbations of the form  $C(x)=\delta C_{0}(x)$ for appropriate $C_{0}(x)$ and
$\delta\in\mathbb{R}$, we find that the rate function for the empirical measure is quadratic in  $\delta\in\mathbb{R}$,  see Theorem
\ref{T:measure2}.

Our second key result is that the information in $I_C(\mu)$ can be used to study specific observable: from the large deviation for the empirical measure
we have a large deviation for principle for observables $f \in \mathcal{C}(E; \mathbb{R})$,
\begin{equation*}
\mathbb{P} \left\{ \frac{1}{t}\int_0^t f(X_s) \, ds  \approx \ell \right\} \asymp e^{- t \tilde{I}_{f,C}(\ell)}
\end{equation*}
and we show that $\tilde{I}_{f,C}(\ell) > \tilde{I}_{f,0}(\ell)$ unless $f$ and $\ell$ satisfy the
non degeneracy condition (in form of a Poisson equation) given in Theorem \ref{T:observable1}, see also Remarks \ref{R:ConditionObservable} and \ref{R:PoissonEquation}.

Moreover, one can deduce information about  asymptotic variances from the large deviations rate function,
since the second derivative of the rate function $\tilde{I}_{f,C}(\ell)$
evaluated at $\ell=\bar{f}$  is inversely proportional to the asymptotic variance of the estimator, denoted by $\sigma_{f,C}^2$.
Based on this relation, we show that the asymptotic variance strictly decreases  $\sigma^2_{f,C} < \sigma^2_{f,0}$, for generic
observables.

The paper is organized as follows. In Section \ref{S:MainResults} we recall some well-known results about large deviations due to
Donsker-Varadhan and G\"artner  and we present our main results. Proofs of statements related to the rate function for
 the empirical measure are in Section  \ref{S:DonskerVaradhanRateFcn}.  In particular, we prove Theorems \ref{T:measure1} and
\ref{T:measure2} by using a  representation of the rate function $I(\mu)$ due to G\"artner \cite{Gartner1977}.
Proofs related to the rate function for a given  observable and the results for variance reduction are in Section
\ref{S:LDPobservable}. In particular, we use the results of Section \ref {S:DonskerVaradhanRateFcn} to deduce the results on
the rate function and asymptotic variance for observables, i.e. Theorems \ref{T:observable1} and \ref{T:VarianceObservable}.
In Section \ref{S:Simulations} we present a few simulation results to illustrate the theoretical findings.

\section{Main results}\label{S:MainResults}

Let us first recall the definition of the large deviations principle for a family of empirical measures $\pi_{t}$.
Let $E$ be a Polish space, i.e., a complete and separable metric space. Denoting by
$\mathcal{P}(E)$ the space of all probability measures on $E$,  we equip $\mathcal{P}(E)$ with the topology of
weak convergence, which makes $\mathcal{P}(E)$  metrizable and a Polish space.

\begin{definition}\label{Def:LDP}
Consider a sequence of random probability measures $\{\pi_{t}\}$. The family $\{\pi_{t}\}$ is said to satisfy a large deviations principle (LDP) with rate function (equivalently action functional) $I:\mathcal{P}(E)\mapsto [0,\infty]$ if the following conditions hold:
\begin{itemize}
\item{For all open sets $O\subset \mathcal{P}(E)$, we have
\[
\liminf_{t\rightarrow\infty}\frac{1}{t}\log \mathbb{P}\left\{\pi_{t}\in O\right\}\geq -\inf_{\mu\in O}I(\mu)
\] }
\item{For all closed sets $F\subset \mathcal{P}(E)$, we have
\[
\limsup_{t\rightarrow\infty}\frac{1}{t}\log \mathbb{P}\left\{\pi_{t}\in F\right\}\leq -\inf_{\mu\in F}I(\mu)
\] }
\item{The level sets $\{\mu: I(\mu)\leq M\}$ are compact in $\mathcal{P}(E)$ for all $M<\infty$.}
\end{itemize}
\end{definition}

If the random measures $\pi_{t}$ are the empirical measures of an ergodic Markov process $X_{t}$ (see \eqref{emp-measure})
with invariant distribution $\bar{\pi}$ then $I(\mu)$ is a nonnegative convex function with $I(\bar{\pi})=0$ and
thus $I(\mu)$ controls the rate at which the random measure $\pi_t$ concentrates to $\bar{\pi}$.

For convenience we will assume that the  diffusion process $X_t$  which solves the SDE \eqref{irreversibleS}
takes values in a compact space and that the vector fields are sufficiently smooth. We fully expect, though, our result to still
hold in  $\mathbb{R}^d$ under suitable confining assumptions on the potential $U$ to ensure a large deviation principle.
Throughout the rest of  the paper we assume that

\bigskip
\noindent{\bf{(H)}}  The state space $E$ is a connected, compact, d-dimensional smooth Riemann manifold without boundary,
and there exists an $\alpha\in(0,1)$ such that
the potential $U \in {\mathcal C}^{(2+\alpha)}(E )$ and the vector field $C \in {\mathcal C}^{(1+\alpha)}(E)$.
Moreover,  we assume that ${\rm div}(C e^{-2U}) =0$ so that the measure ${\bar \pi}$ is invariant.

\bigskip
From the work of G\"artner and Donsker-Vardhan, \cite{Gartner1977,DonskerVaradhan1975},  under condition
${\bf (H)}$, the empirical measures  $\pi_t$ satisfy a large deviation principle which is uniform in the initial condition, i.e.
the rate function is independent of the distribution of $X_0\sim \mu_0$. Let us denote by  $\mathcal{L}$ the infinitesimal generator of the Markov process $X_{t}$ and by $\mathcal{D}$ its domain of definition.
The rate function $I(\mu)$ (usually referred to as the
Donsker-Vardhan functional) takes the form
\begin{equation*}
I(\mu)=-\inf_{u\in\{u\in\mathcal{D},u>0\}}\int_{E}\frac{\mathcal{L}u}{u}d\mu.
\end{equation*}
 An alternative formula for $I(\mu)$, more
useful in the context of this paper,  is given in terms of the Legendre transform
\begin{equation*}
I(\mu)= \sup_{f \in \mathcal{C}(E)} \left\{  \int f \, d\mu - \lambda(f) \right\}\,,
\end{equation*}
where $\lambda(f)$ is the maximal eigenvalue of the Feyman-Kac semigroup
$T_t^f h(x) =  \mathbb{E}_x [ e^{\int_0^t f(X_s) ds} h(X_t)]$ acting on the Banach space $\mathcal{C}(E; \mathbb{R})$.
As shown in  \cite{Gartner1977} for nice $\mu$ this formula can be used to derive a useful, more explicit, formula for $I(\mu)$ which
will be central in our analysis (see Theorem \ref{Th:Gartner} below).

In the sequel and in order to emphasize the dependence on $C$ of the rate function we will use the notation
$I_C(\mu)$. Our first two results show that adding an irreversible drift $C$ increases the
Donsker-Varadhan rate function pointwise.

\begin{theorem}\label{T:measure1} Assume that $C \not=0$ is as in Assumption ${\bf (H)}$. For any $\mu \in \mathcal{P}(E)$ we have $I_C(\mu) \ge
I_0(\mu)$. Let  $\mu(dx) =p(x)dx$ be a probability measure with positive density $p \in \mathcal{C}^{(2 + \alpha)}(E)$ for some $\alpha >0$ and $\mu \not = \bar{\pi}$.
Then, we have
\begin{equation*}
I_C(\mu) = I_0(\mu) + \frac{1}{2}\int_{E}\left|\nabla \psi_{C}(x)-\nabla U(x)\right|^{2}d\mu(x) \,.
\end{equation*}
where $\psi_C$ is the unique solution (up to a constant) of the elliptic equation
\begin{equation*}
{\textrm div}\left[p\left(-\nabla U+C+\nabla \psi_{C}\right)\right]=0.
\end{equation*}

Moreover, we have $I_C(\mu) = I_0(\mu)$ if and only if the positive density  $p(x)$ satisfies $\text{div}\left(p(x)C(x)\right)=0$.
Equivalently such $p$ have the form $p(x) = e^{2 G(x)}$ where $G$ is such that $G+U$ is an invariant  for the
vector field $C$ (i.e., $C \nabla (G+U) =0$).
\end{theorem}

To obtain a slightly more quantitative result let us consider a one-parameter family $C(x) = \delta C_0(X)$ where
$\delta \in \mathbb{R}$ and $C_0$. We show that for any fixed measure $\mu$ the functional $I_{\delta C_0}(\mu)$
is quadratic in $\delta\in\mathbb{R}$.

\begin{theorem}\label{T:measure2} Assume that $C=\delta C_0 \not=0$ is  as in Assumption ${\bf (H)}$ and consider the measure $\mu(dx) =p(x)dx$
with positive density $p \in \mathcal{C}^{(2 + \alpha)}(E)$ for some $\alpha >0$.  Then we have
\begin{equation*}
I_{\delta C_0}(\mu) = I_0(\mu) + \delta^2 K(\mu)  \,,
\end{equation*}
where the functional $K(\mu)$ is strictly positive if and only if $\text{div}\left(p(x)C_0(x)\right)\not=0$. Moreover, the functional $K(\mu)$ takes the explicit form
\begin{equation*}
K(\mu)=\frac{1}{2}\int_{E}\left|\nabla \xi(x)\right|^{2}d\mu(x) \,,
\end{equation*}
where $\xi$ is the unique solution (up to a constant) of the elliptic equation
\begin{equation*}
{\textrm div}\left[p\left(C_{0}+\nabla \xi\right)\right]=0\,.
\end{equation*}
\end{theorem}

For $f \in {\mathcal C}(E)$ the contraction principle implies that the ergodic average $\frac{1}{t} \int_0^t f(X_s) ds$
satisfies a large deviation  principle with the rate function
\begin{equation*}
\tilde{I}_{f,C}(\ell)=\inf_{\mu\in\mathcal{P}(E)}\left\{I_C(\mu): \left<f,\mu\right>=\ell\right\} \,.
\end{equation*}
Note that $\tilde{I}_{f,C}(\ell)$ can also be expressed in terms of a Legendre transform
\begin{equation*}
\tilde{I}_{f,C}(\ell) \,=\, \sup_{\beta \in \mathbb{R}} \left\{  \beta \ell - \lambda(\beta f) \right\}  \,,
\end{equation*}
where
\[
 \lambda(\beta f) =\lim_{t\rightarrow\infty}\frac{1}{t} \log \mathbb{E} \left[ e^{\int_0^t \beta f(X_s) ds}\right].
\]
The eigenvalue $\lambda(\beta f)$ is a smooth strictly convex function of $\beta$ so that if $\ell$ belongs to the
range of $f$ we have
\begin{equation*}
\tilde{I}_{f,C}(\ell) \,=\,  \widehat{\beta} \ell - \lambda(\widehat{\beta} f) \,, \quad {\rm with~} \widehat{\beta} {\rm ~given~ by~} \ell =
\frac{d}{d \beta}  \lambda( \widehat{\beta} f) \,.
\end{equation*}
In fact, if $f \in \mathcal{C}^{(\alpha)}(E)$, then by Proposition \ref{P:ExistenceMinimizer}  there is $\mu^{*}_{C}(dx)=p_{C}(x)dx$, with $p_{C}(x)>0$ and
$p_{C} \in \mathcal{C}^{(2+\alpha)}(E)$ such that $\tilde{I}_{f,C}(\ell)=I_{C}(\mu^{*}_{C})$. Then, Theorem \ref{T:measure1} and Proposition
\ref{P:ExistenceMinimizer} give Theorem \ref{T:observable1}. Theorem \ref{T:observable1} shows that the rate function for observables
increases pointwise under a non-degeneracy condition.

\begin{theorem}\label{T:observable1} Assume that $C\not=0$ is as in Assumption ${\bf (H)}$. Consider
$f \in \mathcal{C}^{(\alpha)}(E)$ and $\ell \in ( \min_x f(x), \max_x f(x))$ with $\ell \not= \int f d \bar{\pi}$.
Then we have
\begin{equation*}
{\tilde I}_{f,C} (\ell) \ge {\tilde I}_{f,0}(\ell) \,.
\end{equation*}

Moreover if there exists $\ell_0$ such that for the vector field $C$, ${\tilde I}_{f,C} (\ell_0) = {\tilde I}_{f,0}(\ell_0)$ then we must have
\begin{equation}
\widehat{\beta}(\ell_0) f \,=\,  \frac{1}{2} \Delta(G+U) + \frac{1}{2} |\nabla G|^2 - \frac{1}{2} |\nabla U|^2 \,,\label{Eq:ConditionObservable}
\end{equation}
where $G$ is such that $G+U$ is invariant under the particular vector field $C$.
\end{theorem}

The following remarks are of interest.
\begin{remark}\label{R:ConditionObservable}
Letting $\mathcal{L}_{0}$ denote the infinitesimal generator of the reversible process $X_{t}$ defined in (\ref{reversibleS}), we get that (\ref{Eq:ConditionObservable}) can be rewritten as a nonlinear Poisson equation of the form
\begin{equation*}
\widehat{\beta}(\ell_0) f \,=\, \mathcal{H}(G+U)\,,
\end{equation*}
where
\[\mathcal{H}(G+U)=e^{-(G+U)}\mathcal{L}_{0}e^{G+U}=\frac{1}{2} \Delta(G+U) + \frac{1}{2} |\nabla G|^2 - \frac{1}{2} |\nabla U|^2.
\]

Recalling Theorem \ref{T:measure1} (see the proof of Theorem \ref{T:observable1}), an alternative condition that gives ${\tilde I}_{f,C} (\ell_0) = {\tilde I}_{f,0}(\ell_0)$  
is as follows.
By Proposition \ref{P:ExistenceMinimizer} there is $\mu^{*}_{C}(dx;\ell_{0})=p_{C}(x;\ell_{0})dx$, with $p_{C}>0$ and
$p_{C} \in \mathcal{C}^{(2+\alpha)}(E)$ such that $\tilde{I}_{f,C}(\ell)=I_{C}(\mu^{*}_{C}(\cdot;\ell_{0}))$. Then, the condition
$\text{div}(p_{C}(x;\ell_{0})C(x))=0$,  implies that ${\tilde I}_{f,C} (\ell_0) = {\tilde I}_{f,0}(\ell_0)$.
\end{remark}

\begin{remark}\label{R:PoissonEquation}
In is interesting to note here that the Poisson equation (\ref{Eq:ConditionObservable}) is reminiscent of Poisson equations that have appeared in the literature in the analysis of MCMC algorithms, see for example Chapter 17 of \cite{MeynTweedie}. In this paper, we see that the particular  Poisson equation can be also used to characterize when irreversible perturbations do actually strictly improve convergence to equilibrium.
\end{remark}

A standard measure of efficiency of a sampling method for an observable $f$ is to use the asymptotic variance.
Under our assumptions the central limit theorem holds for the ergodic average $f_t$ and we have
\begin{equation}
t^{1/2} \left( \frac{1}{t} \int_0^t f(X_s) ds - \int f d \bar{\pi} \right)  \Rightarrow N(0,\sigma_{f,C}^{2})\label{Eq:CLT}
\end{equation}
and the asymptotic variance $\sigma_{f,C}^2$ is given in terms of the integrated autocorrelation function, see e.g., Proposition IV.1.3 in \cite{AsmussenGlynn2007},
\[
\sigma_{f,C}^2 \,=\, 2\int_0^\infty  \mathbb{E}_{\bar{\pi}}\left[\left(f(X_{0})-\bar{f}\right)\left(f(X_{t})-\bar{f}\right)\right] \, dt.
\]

This is a convenient quantity from a practical point of view since there exists easily implementable
estimators for $\sigma^{2}_{f,C}$.  On the other hand the asymptotic variance $\sigma^{2}_{f,C}$ is related to the curvature of the rate function
$I_{f,C} (\ell)$ around the mean $\ell=\bar{f}$ (e.g., see \cite{Hollander2000}): we have
\[
{\tilde I}_{f,C}''( \bar{f}) = \frac{1}{2 \sigma_{f,C}^2}  \,.
\]

From Theorem \ref{T:observable1} it follows immediately that $\sigma^{2}_{f,C} \leq \sigma^{2}_{f,0}$
but in fact  the addition of an appropriate irreversible drift strictly decreases the asymptotic variance.

\begin{theorem}\label{T:VarianceObservable} Assume that $C \not=0$ is a vector field as in assumption $\bf{(H)}$ and let $f \in \mathcal{C}^{(\alpha)}(E)$ such that for some $\epsilon >0$
and $\ell \in ( {\bar f} - \epsilon, {\bar f}+\epsilon) \setminus \left\{\bar{f}\right\}$ we have ${\tilde I}_{f,C} (\ell) > {\tilde I}_{f,0}(\ell)$.
Then we have
\begin{equation*}
\sigma_{f,C}^{2} < \sigma_{f,0}^{2}.
\end{equation*}
\end{theorem}

\begin{remark}
An examination of the proof of Theorem \ref{T:VarianceObservable} shows that a less restrictive condition is needed for the strict decrease in variance to hold. In particular, it is enough to assume that
\begin{equation*}
\text{div}\left(\frac{\partial p_{C}(x) }{\partial \ell}\Big |_{\ell=\bar{f}} C(x)\right)\neq 0
\end{equation*}
where $p_{C}(x)=p_{C}(x;\ell)$ is the strictly positive invariant density of $\mu^{*}_{C}(dx)=\mu^{*}_{C}(dx;\ell)$ such that ${\tilde I}_{f,C} (\ell)=I_{C}(\mu^{*}_{C})$.
\end{remark}

Let us conclude this section with an example demonstrating that adding irreversibility in the dynamics does not always result in a increase of the spectral gap, even though the variance of the estimator decreases.
The key point is that the imaginary part of  complex eigenvalues of the generator for irreversible processes creates oscillations in the autocorrelation function which can dramatically reduce the value of
its integral. A related discussion regarding comparison of convergence criteria can be also found in \cite{DupuisDoll1}. Related computations for the asymptotic behavior of the mean-square displacement of tracers can be found in \cite{MajdaKramer}. The purpose of this example is to demonstrate that spectral gap as a criterium of convergence may not be tight enough to assess improvement in performance when breaking irreversibility. On the other hand, the large deviations rate function and the asymptotic variance both reflect the improved convergence properties due to the irreversible perturbation.

\begin{example}\label{Ex:Example1}
Let us consider the family of diffusions
\[
dX_{t}=  \delta dt+dW_{t}
\]
on the circle $S^1$ with generator
\[
\mathcal{L}_{\delta} \,=\, \Delta + \delta \nabla
\]

For any  $\delta\in\mathbb{R}$ the Lebesgue measure  on $S^1$ is invariant, but $\mathcal{L}_{\delta}$ is self-adjoint  on
$L^2(dx)$ and thus $X_t$ is reversible if  and only if $\delta=0$.   A simple computation (using for example Lemma \ref{Th:Gartner1}) shows that for a measure $\mu(dx)$ that has positive and sufficiently smooth density $p(x)$ we have
\[
I(\mu)=\frac{1}{8}\int_{S^{1}}\left|\frac{ p'(x)}{p(x)}\right|^{2}p(x)dx+\delta^{2}\frac{1}{2}\left[1-
\frac{1}{\int_{S^{1}}\frac{1}{p(x)}dx}\right],
\]
and in this case $I(\mu)$ strictly increases unless $\mu(dx)= dx$. The eigenvalues of  $\mathcal{L}_\delta$
are $ \lambda_n =  -n^2 + in\delta, \quad  n \in \mathbb{Z}$ with eigenfunction $e^{inx}$ and thus the spectral gap is $-1$  for any
$\delta\in\mathbb{R}$.  However for any real-valued function $f$ the asymptotic variance decreases:  for $f$ with $\int_{S^1}f dx =0$
with Fourier coefficients $c_n$ we have
\[
\sigma^{2}_{f}(\delta)= \int_0^\infty   \langle e^{t {\mathcal L}} f(x)\,,  f(x) \rangle_{L^2(dx)}  \, dt\,=\, \sum_{n \in \mathbb{Z}, n\not=0} \frac{|c_n|^2} {n^2 +in\delta} \,=\, \sum_{n=1}^\infty \frac{ 2|c_n|^2}{n^2 + \delta^2}.
\]
In this example, even though the spectral gap does not increase at all, the variance not only decreases, but it can be made as
small as we want by increasing $\delta^{2}$.  The latter is in agreement with  both Theorem \ref{T:measure2} and Theorem
\ref{T:VarianceObservable} and illustrates how irreversibility improves sampling.
\end{example}

\section{The Donsker-Vardhan functional}\label{S:DonskerVaradhanRateFcn}

A standard trick in the theory of large deviations, when computing the probability of an unlikely event,  is to
perform a change of measure to make the unlikely event typical. In the context of SDE's, this takes of
the form of changing the drift of the SDE's itself.  This is the idea behind the proof of the following result
due to G\"artner,  \cite{Gartner1977}.

\begin{theorem}[Theorem 3.2 in \cite{Gartner1977}] \label{Th:Gartner} Consider the SDE
\[
dX_t =  b(X_t) + dW_t \,
\]
on $E$ with $b \in {\mathcal C}^{(1+\alpha)}(E)$ and with generator
\[
\mathcal{L} = \Delta + b \nabla \,.
\]

Let $\mu \in \mathcal{P}(E)$, where  $\mu(dx) =p(x)dx$ is a measure with positive density $p \in \mathcal{C}^{(2 + \alpha)}(E)$ for some $\alpha >0$.
The Donsker-Vardhan rate function $I(\mu)$ takes the form
\begin{equation}\label{Eq:GartnerFormula1}
I(\mu)=\frac{1}{2}\int_{E}\left|\nabla \phi(x)\right|^{2}d\mu(x) \,,
\end{equation}
where $\phi$ is the unique (up to constant) solution of the equation
\begin{equation}\label{Eq:GartnerFormula1Constraint}
\Delta\phi+\frac{1}{p}\left(\nabla p,\nabla\phi\right)=\frac{1}{p}\mathcal{L}^{*}p \,,
\end{equation}
and $\mathcal{L}^*= \Delta - \nabla b$ is the formal adjoint of $L$ in $L^2(dx)$.
\end{theorem}

In the special case where $b=-\nabla U$ is a gradient, then up to an additive constant $\phi(x)=\frac{1}{2}\log p(x)+ U(X)$,
and we get
\begin{equation}\label{explicitreversible}
I(\mu)=\frac{1}{2}\int_{E}\left|\frac{1}{2}\frac{\nabla p(x)}{p(x)}+\nabla U(x)\right|^{2}d\mu(x)
\end{equation}
which is the usual explicit formula for the rate function in the reversible case.

It will be useful to rewrite $I(\mu)$ in a  different form.
\begin{lemma}\label{Th:Gartner1} Under the conditions of Theorem \ref{Th:Gartner}, we have
\begin{equation*}
I(\mu)=\frac{1}{8}\int_{E}\left|\frac{\nabla p(x)}{p(x)}\right|^{2}d\mu(x)+\frac{1}{2}\int_{E}\left|\nabla \psi(x)\right|^{2}d\mu(x)-\frac{1}{2}\int_{E}\frac{b\nabla p}{p}d\mu(x) \,,
\end{equation*}
where $\psi$ is the unique (up to constant) solution of the elliptic equation
\begin{equation*}
\textrm{div}\left[p\left(b+\nabla \psi\right)\right]=0 \,.
\end{equation*}
\end{lemma}

\proof  Motivated by the solution in gradient case, let us write $\phi(x)=\frac{1}{2}\log p(x)+\psi(x)$. By plugging  $\phi(x)=\frac{1}{2}\log p(x)+\psi(x)$ in (\ref{Eq:GartnerFormula1}), we get
\begin{align}
I(\mu)&=\frac{1}{2}\int_{E}\left|\frac{1}{2}\frac{\nabla p(x)}{p(x)}+\nabla \psi(x)\right|^{2}d\mu(x)\nonumber\\
&=\frac{1}{8}\int_{E}\left|\frac{\nabla p(x)}{p(x)}\right|^{2}d\mu(x)+\frac{1}{2}\int_{E}\left|\nabla \psi(x)\right|^{2}d\mu(x)+\frac{1}{2}\int_{E}\frac{\nabla\psi \nabla p}{p}d\mu(x)\nonumber\\
&=\left[\frac{1}{8}\int_{E}\left|\frac{\nabla p(x)}{p(x)}\right|^{2}d\mu(x)+\frac{1}{2}\int_{E}\left|\nabla \psi(x)\right|^{2}d\mu(x)-\frac{1}{2}\int_{E}\frac{b\nabla p}{p}d\mu(x)\right]+\frac{1}{2}\int_{E}\left[\left(b+\nabla\psi\right)\nabla p\right] dx\nonumber\\
&= I(\mu,1)+I(\mu,2)\nonumber \,.
\end{align}

We claim that $I(\mu,2)=0$. Indeed, using $\phi(x)=\frac{1}{2}\log p(x)+\psi(x)$, the constraint  (\ref{Eq:GartnerFormula1Constraint}) gives the following chain of equalities
\begin{align}
\Delta\phi+\frac{1}{p}\left(\nabla p,\nabla\phi\right)&=\frac{1}{p}\mathcal{L}^{*}p\Rightarrow\nonumber\\
\frac{\Delta p}{2p}-\frac{|\nabla p|^{2}}{2p^{2}}+\Delta \psi+ \frac{|\nabla p|^{2}}{2p^{2}}+\frac{1}{p}\left(\nabla p,\nabla\psi\right)&= \frac{\Delta p}{2p}-\frac{1}{p}\textrm{div}(bp)\Rightarrow\nonumber\\
\Delta \psi+\frac{1}{p}\left(\nabla p,\nabla\psi\right)&= -\frac{1}{p}\textrm{div}(bp)\Rightarrow\nonumber\\
p\Delta \psi+\left(\nabla p,\nabla\psi\right)&+\textrm{div}(bp)=0\Rightarrow\nonumber\\
\nabla \cdot \left[p\left(b+\nabla \psi\right)\right]&=0\nonumber \,.
\end{align}

The weak formulation of the latter statement reads as follows
\[
\int_{E}\left(b(x)+\nabla \psi(x)\right)\nabla g(x)p(x)dx=0,\quad \forall g\in\mathcal{C}^{1}(E) \,.
\]

Choosing $g=\log p$, we  obtain
\[
\int_{E}\left(b(x)+\nabla \psi(x)\right)\nabla p(x)dx=0,
\]
which is precisely the statement $I(\mu,2)=0$. So we have indeed proven the claim. \qed

\bigskip

With the representation of $I_C(\mu)$ we can now prove  Theorem \ref{T:measure1}.  \medskip

\noindent
{\em Proof of Theorem \ref{T:measure1}}: Since $b(x)=-\nabla U(x)+C(x)$,  using Lemma \ref{Th:Gartner1}, $I_{C}(\mu)$
becomes
\begin{eqnarray}
I_{C}(\mu)&=&\frac{1}{8}\int_{E}\left|\frac{\nabla p(x)}{p(x)}\right|^{2}d\mu(x)+\frac{1}{2}\int_{E}\left|\nabla \psi_{C}(x)\right|^{2}d\mu(x)
\nonumber \\
&& +\frac{1}{2}\int_{E}\frac{\nabla U(x)\nabla p(x)}{p(x)}d\mu(x)-\frac{1}{2}\int_{E}\frac{C(x)\nabla p(x)}{p(x)}d\mu(x)\label{Eq:GartnerFormula2C} \,,
\end{eqnarray}
where $\psi_{C}$ is the unique (up to constant) solution of the equation
\begin{equation*}
\textrm{div}\left[p\left(-\nabla U+C+\nabla \psi_{C}\right)\right]=0.
\end{equation*}

The proof of Lemma \ref{Th:Gartner1} shows that $\psi_{C}(x)=\phi(x)-\frac{1}{2}\log p(x)$ where $\phi$ is the unique solution (up to constants) of the equation
(\ref{Eq:GartnerFormula1Constraint}) with $\mathcal{L}=\mathcal{L}_{0}+C\nabla$.

Using the explicit formula \eqref{explicitreversible} for the reversible case we obtain for the difference
$J_C(\mu) = I_C(\mu)- I_0(\mu)$
\begin{align*}
J_C(\mu) &= I_{C}(\mu) - I_0(\mu)
\,=\, \frac{1}{2}\int_{E}\left[\left|\nabla \psi_{C}(x)\right|^{2}-\left|\nabla U(x)\right|^{2}\right]d\mu(x)-\frac{1}{2}\int_{E}\frac{C(x)\nabla p(x)}{p(x)}d\mu(x).\label{Eq:GartnerFormula3C}
\end{align*}

The condition
$\textrm{div}\left(C(x)e^{-2U(x)}\right)=0$ can be rewritten as
\begin{equation*}
\textrm{div}C(x)=2 C(x)\nabla U(x) \,.
\end{equation*}

Integration by parts gives for the last term in $J_{C}(\mu)$
\begin{align*}
\int_{E}\frac{C(x)\nabla p(x)}{p(x)}d\mu(x) &=\int_{E} C(x)\nabla p(x)dx=-\int_{E} \textrm{div}C(x) p(x)dx=-\int_{E} \textrm{div}C(x) d\mu(x)\nonumber\\
&=-\int_{E} 2 C(x)\nabla U(x) d\mu(x).
\end{align*}

Hence, we obtain
\begin{equation*}
J_{C}(\mu)=\frac{1}{2}\int_{E}\left[\left|\nabla \psi_{C}(x)\right|^{2}-\left|\nabla U(x)\right|^{2}+2 C(x)\nabla U(x)\right]d\mu(x)\label{Eq:RepDifference1} \,.
\end{equation*}

Using the constraint in its weak form
\begin{equation}
 \int_{E}\left[\nabla \psi_{C}(x)-\nabla U(x)+C(x)\right]\nabla g(x)d\mu(x)=0, \quad \textrm{for every }g\in\mathcal{C}^{1}(E)\label{Eq:WeakForm1}
\end{equation}
we can pick freely $g\in\mathcal{C}^{1}(E)$. If we first choose $g=\psi_{C}+U$, then, (\ref{Eq:WeakForm1}) gives
\begin{equation*}
 \int_{E}\left[\left|\nabla \psi_{C}(x)\right|^{2}-\left|\nabla U(x)\right|^{2}\right]d\mu(x)=-\int_{E}C(x)\left(\nabla\psi_{C}(x)+\nabla U(x)\right)d\mu(x)
\end{equation*}
and thus
\begin{equation}
J_{C}(\mu)=\frac{1}{2}\int_{E} C(x)\left(\nabla U(x)-\nabla\psi_{C}(x)\right)d\mu(x).\label{Eq:RepDifference2}
\end{equation}

Choosing $g=\psi_{C}-U$ and we get from (\ref{Eq:WeakForm1})
\begin{equation*}
 \int_{E}\left|\nabla \psi_{C}(x)-\nabla U(x)\right|^{2}d\mu(x)=\int_{E}C(x)\left(\nabla U(x)-\nabla\psi_{C}(x)\right)d\mu(x)\,.
\end{equation*}
Plugging this in (\ref{Eq:RepDifference2}) we obtain
\begin{equation*}
J_{C}(\mu)=\frac{1}{2}\int_{E}\left|\nabla \psi_{C}(x)-\nabla U(x)\right|^{2}d\mu(x) \,.\label{Eq:RepDifference3}
\end{equation*}

Clearly $J_{C}(\mu)\geq 0$. If $\mu$ possesses a strictly positive density, it is clear that
$J_{C}(\mu)= 0$ if and only if $\text{div}\left(pC\right)= 0$. In other words, $J_{C}(\mu)> 0$ if and only if $\text{div}\left(pC\right)\neq 0$.

Finally  let us write the positive density as $p(x)= e^{ 2 G(x)}$, since we have $\text{div} (C e^{-2 U})=0$ and $\text{div}( C
e^{2 G})=0$ we have
\[
{\textrm{ div}} C = -2 C \nabla U = 2 C \nabla G
\]
and thus $C \nabla (G+U) =0$, i.e. $(G+U)$ is a conserved quantity under the flow $\frac{dx}{dt} = C(x)$.
\qed

\bigskip

We now consider the one-parameter family $C(x) = \delta C_0(x)$ and prove Theorem \ref{T:measure2}.

\bigskip
\noindent
{\em Proof of Theorem \ref{T:measure2}:}   For notational convenience let us write $J_\delta(\mu)$ instead of $J_{\delta C_0}(\mu)$ and
let us set $\varphi_{\delta}(x)=\psi_{\delta C_0}(x)- U(x)$. From Theorem  \ref{T:measure1} we have
\begin{equation}\label{Eq:GartnerJa}
J_{\delta}(\mu)=\frac{1}{2}\int_{E}\left|\nabla \varphi_{\delta}(x)\right|^{2}d\mu(x)
\end{equation}
where $\varphi_{\delta}$ is the unique (up to constant) solution of the equation
\begin{equation}
\int_{E}\left(\delta C_0(x)+\nabla \varphi_{\delta}(x)\right)\nabla g(x)\mu(dx)=0,\quad \forall g\in\mathcal{C}^{1}(E).\label{Eq:GartnerFormula2ConstraintC2}
\end{equation}

Let us define $\xi_{\delta}(x)=\delta^{-1} \varphi_{\delta}(x)$. Then,
\begin{equation*}
J_{\delta}(\mu)=\delta^{2}\frac{1}{2}\int_{E}\left|\nabla \xi_{\delta}(x)\right|^{2}d\mu(x)
\end{equation*}
and because $\delta\neq0$, $\xi_{\delta}$ is the unique (up to constant) solution of the equation
\begin{equation*}
\int_{E}\left(C_0(x)+\nabla \xi_{\delta}(x)\right)\nabla g(x)\mu(dx)=0,\quad \forall g\in\mathcal{C}^{1}(E).
\end{equation*}

The last equation makes it clear that, modulo an additive constant, $\xi_{\delta}(x)$ is in fact independent of $\delta$.
 Thus, there exists a functional $K(\mu)\geq 0$ such that
\begin{equation*}
J_{\delta}(\mu)=\delta^{2}K(\mu) \,.
\end{equation*}
Clearly, if $\mu(dx)=p(x)dx$ with $\textrm{div}(pC_{0})=0$ then $K(\mu)=0$, otherwise $K(\mu)> 0$.
\qed

\section{Large deviation for observables and the asymptotic variance}\label{S:LDPobservable}

Let us consider a function $f\in \mathcal{C}(E)$ with mean $\bar{f}=\int_{E}f(x)d\bar{\pi}(x)$. Let us set
\begin{equation*}
f_{t}=\left<f,\pi_{t}\right>=\int_{E}f(x)d\pi_{t}(x) =\frac{1}{t}\int_{0}^{t}f(X_{s})ds \,.
\end{equation*}

By the contraction principle $f_{t}$ satisfies a large deviation principle with action functional
given by
\begin{equation}\label{Eq:ContractionPrinciple}
\tilde{I}_{f,C}(\ell)=\inf_{\mu\in\mathcal{P}(E)}\left\{I_C(\mu): \left<f,\mu\right>=\ell\right\} \,,
\end{equation}
where $\ell\in\mathbb{R}$ and $I_C(\mu)$ is the Donsker-Vardhan action functional for the empirical measure $\pi_{t}$.

In Subsection \ref{SS:LDPobservable} we prove Theorem \ref{T:observable1}, whereas in Subsection \ref{SS:AsymptoticVariance} we prove Theorem \ref{T:VarianceObservable}.

\subsection{Large deviation for observables}\label{SS:LDPobservable}
Theorem \ref{T:observable1} is a fairly immediate consequence of Theorem \ref{T:measure1} and Proposition
\ref{P:ExistenceMinimizer}.

\begin{proposition}\label{P:ExistenceMinimizer}
Let $f \in \mathcal{C}^{(\alpha)}(E)$, and $\ell \in ( \min_x f(x), \max_x f(x))$. Then there exists
$\mu^*(dx) = p(x)dx$  with $p(x) > 0$ and $p(x)  \in \mathcal{C}^{(2+\alpha)}(E)$ such that
\begin{equation*}
\tilde{I}_{f,C}(\ell) = I_{C}(\mu^*) \,.
\end{equation*}
\end{proposition}

\proof  As discussed in G\"artner \cite{Gartner1977}, the semigroup
$T_t h(x) = \mathbb{E}_x \left[ h(X_t) \right]$  is strong-Feller and the strong-Feller property is inherited by
the Feynman-Kac semigroup
\begin{equation*}
T^f_t h(x) \,=\, \mathbb{E}_x \left[ e^{\int_{0}^{t} f(X_{s})ds} h(X_t) \right] \,,
\end{equation*}
if $f \in \mathcal{C}(E)$.  Moreover the semigroups $T^f_t $ are quasi-compact on the Banach space
$\mathcal{C}(E)$ and by a Perron-Frobenius argument the semigroup $T^f_t$ has a dominant simple
positive eigenvalue $e^{\lambda(f)t}$ with a corresponding  strictly positive eigenvector $u(f)= e^{\phi(f)}$. We write $\lambda(f)$ and $u(f)$ instead of $\lambda,u$ in order to
emphasize their dependence on the observable $f$.

For any $f,g \in \mathcal{C}(E)$, $T_t^{f+\gamma g}$ is a bounded perturbation of $T^f_t$.  By analytic perturbation theory (see for example Chapter VIII of \cite{Kato1969})
and the simplicity of the eigenvalue $\lambda(f)$ this implies that the maps $\gamma \mapsto \lambda(f + \gamma g)$ and
$\gamma \mapsto u(f + \gamma g)$ are real-analytic functions.
If we require, in addition, that $f \in \mathcal{C}^{(\alpha)}(E)$, then the bounded linear operator $(\mathcal{L}_C + f)$  that maps $\mathcal{C}^{(2+\alpha)}(E)$ to
$\mathcal{C}^{(\alpha)}(E)$ is invertible with compact inverse. Hence, the relation
\begin{equation*}
(\mathcal{L}_C + f)  u(f) = \lambda(f) u(f)  \,.
\end{equation*}
implies that $\lambda(f)=\lim_{t\rightarrow\infty} \frac{1}{t}\log \mathbb{E}\left[e^{\int_{0}^{t} f(X_{s})ds}\right]$ is a simple eigenvalue of the operator $(\mathcal{L}_C + f)$ in $\mathcal{C}^{(\alpha)}(E)$ and that the solution $u(f)$ is in $\mathcal{C}^{(2 + \alpha)}(E)$ (see \cite{DouglisNirenberg1955}). This implies
\begin{equation*}
\nabla \phi(f) = \nabla \log u(f) \in \mathcal{C}^{(1+\alpha)}(E) \,.
\end{equation*}

The rate function $I_C(\mu)$ can be written as
\begin{equation*}\label{Eq:RateLegendre1}
I_C(\mu) \,=\, \sup_{f \in \mathcal{C}^{(\alpha)}(E)} \left\{ \mu(f) - \lambda(f) \right\} \,.
\end{equation*}
If we pick $\mu(dx) = p(x) dx$ with $p(x)>0$ and $p \in C^{(2 + \alpha)}(E)$ then
it is shown in \cite{Gartner1977} that the supremum is attained when $f$ is chosen
such that $\mu$ is the invariant measure for the SDE with infinitesimal generator
\begin{equation*}
\mathcal{L}_C + \nabla \phi(f) \nabla = \mathcal{L}_{C + \nabla \phi} \,.
\end{equation*}

Turning now to the rate function for observables we note first that if $ \ell \in (\min_x f(x), \max_x f(x))$ then $I_{f,C}(\ell)$
is finite. Indeed simply pick any measure $\mu$ with a $\mathcal{C}^{(2+\alpha)}(E)$ strictly positive density such that
$\int f d\mu= \ell$, then
$I_{f,C}(\ell) \leq I_C(\mu)$ which is finite by Theorem \ref{Th:Gartner}.
Besides the representation \eqref{Eq:ContractionPrinciple}
we can also represent the rate function  $\tilde{I}_{f,C}$ as the Legendre transform of the moment generating function of
$\bar{f}_{t}$
\begin{equation*}
\tilde{I}_{f,C}(\ell)=\sup_{\beta\in\mathbb{R}}\left\{\ell\cdot\beta-\lambda(\beta f)\right\}
\end{equation*}
where
\begin{equation*}
\lambda(\beta f)=\lim_{t\rightarrow\infty} \frac{1}{t}\log \mathbb{E}\left[e^{\int_{0}^{t}\beta f(X_{s})ds}\right]\,.
\end{equation*}

Due to the relation
\begin{equation}\label{Eq:EigenvalueAlpha1}
(\mathcal{L}_C + \beta f) u(\beta f) = \lambda(\beta f) u(\beta f) \,,
\end{equation}
$\lambda(\beta f)$ is a simple eigenvalue of $\mathcal{L}_C + \beta f$ in $\mathcal{C}^{(\alpha)}(E)$ and as mentioned before  $u(\beta f)$ is in $\mathcal{C}^{(2 + \alpha)}(E)$.
We can then compute $\tilde{I}_{f,C}(\ell)$ by calculus and the sup is attained if $\hat{\beta}$ is chosen such that
$\ell=\frac{\partial}{\partial \beta}  \lambda(\hat{\beta}f)$.
With $u(\beta f) =e^{\phi(\beta f)}$,  the  eigenvalue equation (\ref{Eq:EigenvalueAlpha1}) can be equivalently written as
\begin{equation}\label{Eq:EigenvalueAlpha2}
\mathcal{L}_{C}\phi(\beta f) +\frac{1}{2}\left|\nabla \phi(\beta f)\right|^{2}=\lambda(\beta f)-\beta f
\end{equation}
Differentiating \eqref{Eq:EigenvalueAlpha2} with respect to $\beta$ and setting
$\psi(\beta f) =\frac{\partial \phi}{\partial \beta}(\beta f)$ we see that $\psi(\beta f)$ satisfies
the equation
\begin{equation*}
\mathcal{L}_{C}\psi(\beta f)+\left(\nabla \phi(\beta f),\nabla \psi(\beta f) \right)= \frac{d}{d \beta}\lambda(\beta f)- f
\end{equation*}
or equivalently
\begin{equation*}
\mathcal{L}_{C+\nabla\phi(\beta f)}\psi= \frac{d}{d\beta} \lambda(\beta f)- f
\end{equation*}
Thus, the constraint $\left<f,\mu \right>=\ell$, implies that in order to have $\ell=\frac{d}{d\beta}\lambda(\hat{\beta} f)$ for some $
\hat{\beta}$,  $\mu_{\hat{\beta}}$ should be the invariant measure for the process with generator
$\mathcal{L}_{C+\nabla\phi(\hat{\beta} f)}$.
Since $\nabla \phi \in \mathcal{C}^{(1+\alpha)}(E)$ the corresponding invariant measure $\mu_{\hat{\beta}}$ is strictly positive
and has a density $p(x) \in \mathcal{C}^{(2+\alpha)}(E)$.

To conclude the proof of the proposition, by \cite{Gartner1977} we have $I_C(\mu_{\hat{\beta}}) = \mu( \hat{\beta} f) - \lambda ( \hat{\beta} f)$. But since  $\mu(f) = \ell$ this is also equal to $I_{f,C}(\ell)$.
\qed

\bigskip

\noindent
{\em Completion of the proof of Theorem \ref{T:observable1}}:   Let $\ell$ be such that $\ell \not= \int f d\bar{\pi}$.
By Proposition \ref{P:ExistenceMinimizer},  there exists measures  $\mu^{*}_0$ and $\mu^{*}_C$, both with
strictly positive densities $p_{0},p_{C}\in \mathcal{C}^{(2 + \alpha)}(E)$  such that   $\tilde{I}_{f,C}(\ell)=I_{C}(\mu^{*}_{C})$ and
$\tilde{I}_{f,0}(\ell)=I_{0}(\mu^{*}_{0})$.

Let us first assume that $\textrm{div}(p_{C}C)\neq 0$. Since $I_{C}(\mu) > I_0(\mu)$ for any $\mu$ with strictly positive densities
$p\in C^{(2 + \alpha)}$ such that $\textrm{div}(pC)\neq 0$, this implies that $\tilde{I}_{f,0}(\ell) \le \tilde{I}_{f,C}(\ell)$.

By contradiction let us now assume that
\[
\tilde{I}_{f,0}(\ell) = \tilde{I}_{f,C}(\ell)\,.
\]

Let us first assume that $\mu^{*}_0 \not= \mu^{*}_C$. Since $\textrm{div}(p_{C}C)\neq 0$, we have
\[
I_0(\mu^*_0) = I_C( \mu^*_C) > I_0( \mu^*_C).
\]
which  contradicts $\tilde{I}_{f,0}(\ell) =  I_0(\mu^*_0)$.  Now if $\mu^*_0=\mu^*_C$ then we have
\[
I_{0}(\mu^{*}_{C})=I_{0}(\mu^{*}_{0})=I_{C}(\mu^{*}_{C})\,.
\]

However, this contradicts the fact that we always have $I_{C}(\mu^{*}_{C})> I_{0}(\mu^{*}_{C})$ for $\mu^{*}_C(dx)=p_{C}(x)dx$ such that
$\textrm{div}(p_{C}C)\neq 0$. This proves that $\tilde{I}_{f,0}(\ell) < \tilde{I}_{f,C}(\ell)$.

If $\textrm{div}(p_{C}C) = 0$ then with $p = e^{2G}$ we must have $C \nabla (G+U) =0$.  As in the proof of Proposition \ref{P:ExistenceMinimizer},
the density $p_C$ is an invariant measure for the SDE with added drift $\phi_C$, i.e.,
\[
\mathcal{L}^{*}_{C+\nabla\phi_{C}}p_{C}=0.
\]
but since $\textrm{div} ( p_C C )=0$ we have in fact
\[
\mathcal{L}^{*}_{\nabla\phi_{C}}p_{C}=0.
\]
Also  $\mathcal{L}_{\nabla\phi}$ is the generator of a
reversible ergodic Markov process and thus $p_C = e^{ 2 (\phi - U)}$ from which we see that
\[
\phi = G+U \,.
\]
On the other hand $e^\phi$ is the solution of the eigenvalue equation
\[
({ \mathcal L}_C +  \widehat{\beta} f ) e^\phi   = \lambda(f) e^{\phi}.
\]
Since $C \nabla (G+U) =0$, we have that $C \nabla e^{\phi} =  C \nabla e^{G+U}  =0$. Thus, the last display reduces to
\[
({ \mathcal L}_0 +   \widehat{\beta} f ) e^\phi =  \lambda(f) e^{\phi}.
\]

  We also note that changing $f$ into $f+c$
leaves $\phi$ unchanged but changes $\lambda(f)$ to $\lambda(f) +  \widehat{\beta}c$.
So, for some constant $c$, we must have
\[
\widehat{\beta} f = e^{-(G+U)}   { \mathcal L}_0  e^{(G+U)}  + c \,=\,  \frac{1}{2} \Delta (G +U)  + \frac{1}{2} |\nabla G|^2  - |\nabla G|^2 +
c.
\]

\qed

\subsection{Asymptotic variance}\label{SS:AsymptoticVariance}
In this subsection, we prove that adding irreversibility results in reducing the asymptotic variance of the estimator. The existence of the central limit theorem, see (\ref{Eq:CLT}), of the second derivative $\tilde I_f''( \bar{f})$ and of the relation $\sigma_f^2 = \frac{1}{2 {\tilde I}_f''( \bar{f})}$ implies that it is enough to prove that for $C \not=0$ and $f \in \mathcal{C}^{(\alpha)}(E)$
\begin{equation*}
\tilde{I}_{f,C}^{\prime\prime}(\bar{f})-\tilde{I}_{f,0}^{\prime\prime}(\bar{f})>0
\end{equation*}

We recall that by (\ref{Eq:RepDifference3}),
\begin{align*}
J_C(\mu) &= I_C(\mu)- I_0(\mu)=\frac{1}{2}\int_{E}\left|\nabla \psi_{C}(x)-\nabla U(x)\right|^{2}d\mu(x).
\end{align*}

By Proposition \ref{P:ExistenceMinimizer}, it is enough to consider measures that have a strictly positive density in $\mathcal{C}^{(2+\alpha)}(E)$. We start by computing the first and second order G\^{a}teaux directional derivatives of $J_C(\mu)$ for $\mu(dx)=p(x)dx$ with $p(x)  \in \mathcal{C}^{(2+\alpha)}(E)$. For notational convenience we shall often write $J_C(p)$ instead of $J_C(\mu)$. Let $\gamma\in\mathbb{R}$ and let us define
\begin{equation}
\tilde{J}_{C}(\gamma;p,q)=J_{C}(p+\gamma q), \quad \text{ for } p,q\in\mathcal{C}^{(2+\alpha)}(E).
\end{equation}

In Subsubsection \ref{SSS:FirstGateauxDerivative} we compute first order G\^{a}teaux directional derivative, whereas in Subsubsection \ref{SSS:FirstGateauxDerivative} we compute second order G\^{a}teaux directional derivative. Then, in Subsection \ref{SSS:ProofAsymptoticVariance} we put things together proving Theorem \ref{T:VarianceObservable}.
\subsubsection{First order G\^{a}teaux directional derivative}\label{SSS:FirstGateauxDerivative}
Let $p(x),q(x)  \in \mathcal{C}^{(2+\alpha)}(E)$ and notice that
\begin{align*}
&\frac{1}{\gamma}\left[J_{C}(p+\gamma q)-J_{C}(p)\right]=\nonumber\\
&=\frac{1}{\gamma}\left[\frac{1}{2}\int_{E}\left|\nabla \psi^{p+\gamma q}_{C}(x)-\nabla U(x)\right|^{2}\left(p(x)+\gamma q(x)\right)dx-\frac{1}{2}\int_{E}\left|\nabla \psi_{C}^{p}(x)-\nabla U(x)\right|^{2}p(x)dx\right]\nonumber\\
&=\frac{1}{2\gamma}\left[\int_{E}\left(\nabla \psi^{p+\gamma q}_{C}(x)-\nabla \psi^{p}_{C}(x)\right)\left(\nabla \psi^{p+\gamma q}_{C}(x)+\nabla \psi^{p}_{C}(x)-2\nabla U(x)\right)p(x)dx+\right.\nonumber\\
&\qquad\left.+\gamma\int_{E}\left|\nabla \psi_{C}^{p+\gamma q}(x)-\nabla U(x)\right|^{2}q(x)dx\right]\nonumber\\
&=\frac{1}{2}\left[\int_{E}\frac{\nabla \psi^{p+\gamma q}_{C}(x)-\nabla \psi^{p}_{C}(x)}{\gamma}\left(\nabla \psi^{p+\gamma q}_{C}(x)+\nabla \psi^{p}_{C}(x)-2\nabla U(x)\right)p(x)dx+\right.\nonumber\\
&\qquad\left.+\int_{E}\left|\nabla \psi_{C}^{p+\gamma q}(x)-\nabla U(x)\right|^{2}q(x)dx\right]
\end{align*}

For every $g\in\mathcal{C}^{1}(E)$, we notice that $\frac{\nabla \psi^{p+\gamma q}_{C}(x)-\nabla \psi^{p}_{C}(x)}{\gamma}$ satisfies
\begin{align*}
0&=\frac{1}{\gamma}\int_{E}\left[\left(-\nabla U(x)+C(x)+\nabla \psi^{p+\gamma q}_{C}(x)\right)\nabla g(x)(p(x)+\gamma q(x))\right.\nonumber\\
&\qquad\left.-\left(-\nabla U(x)+C(x)+\nabla \psi^{p}_{C}(x)\right)\nabla g(x)p(x)\right]dx\nonumber\\
&=\int_{E}\frac{\nabla \psi^{p+\gamma q}_{C}(x)-\nabla \psi^{p}_{C}(x)}{\gamma}\nabla g(x) p(x)dx+\int_{E}\left(-\nabla U(x)+C(x)+\nabla \psi^{p+\gamma q}_{C}(x)\right)\nabla g(x)q(x)dx
\end{align*}

Since $p,q\in\mathcal{C}^{(2+\alpha)}(E)$, it follows (as in  Section 3 of \cite{Gartner1977}) that there is a $\hat{\psi}^{p,q}_{C}\in \mathcal{C}^{(2+\alpha)}(E)$ such that
\begin{equation*}
\psi^{p+\gamma q}_{C}(x)=\psi^{p}_{C}(x)+\gamma \hat{\psi}^{p,q}_{C}(x)+o_{1}(\gamma).
\end{equation*}
where $\left\Vert o_{1}(\gamma)\right\Vert_{(2+\alpha)}\rightarrow 0$ as $\gamma\rightarrow 0$. Then $\forall g\in\mathcal{C}^{1}(E)$,  $\nabla \hat{\psi}^{p,q}_{C}(x)$ satisfies
\begin{align}
\int_{E}\left[\nabla \hat{\psi}^{p,q}_{C}(x) p(x)+\left(-\nabla U(x)+C(x)+\nabla \psi^{p}_{C}(x)\right)q(x)\right]\nabla g(x)dx=0. \label{Eq:FirstDerivativeConstraint}
\end{align}

Let us then denote
\[
dJ_{C}(p;q)=\lim_{\gamma\downarrow 0}\frac{J_{C}(p+\gamma q)-J_{C}(p)}{\gamma}.
\]

We obtain
\begin{equation*}
dJ_{C}(p;q)=\int_{E}\nabla\hat{\psi}^{p,q}_{C}(x)\left(\nabla \psi^{p}_{C}(x)-\nabla U(x)\right)p(x)dx+\frac{1}{2}\int_{E}\left|\nabla \psi_{C}^{p}(x)-\nabla U(x)\right|^{2}q(x)dx.
\end{equation*}

It is clear that if the measure $\mu$ is the invariant measure, i.e., $\mu(dx)=\bar{\pi}(dx)$, then denoting by $\bar{p}$ the density of $\bar{\pi}(dx)=\bar{p}(x)dx$, we have that
$\nabla \psi^{\bar{p}}_{C}(x)=\nabla U(x)$.  The latter implies that for any direction $q$, we get
\[
dJ_{C}(\bar{p};q)=0,
\]
which is of course expected to be true.
\subsubsection{Second order G\^{a}teaux directional derivative}\label{SSS:SecondGateauxDerivative}

Next we compute the second order G\^{a}teaux directional derivative. For $p(x),q(x), h(x)  \in \mathcal{C}^{(2+\alpha)}(E)$, we get
\begin{align*}
&\frac{1}{\gamma}\left[dJ_{C}(p+\gamma h;q)-dJ_{C}(p;q)\right]=\nonumber\\
&=\frac{1}{\gamma}\int_{E}\nabla\hat{\psi}^{p+\gamma h,q}_{C}(x)\left(\nabla \psi^{p+\gamma h}_{C}(x)-\nabla U(x)\right)(p(x)+\gamma h(x))dx\nonumber\\
&\qquad-\frac{1}{\gamma}\int_{E}\nabla\hat{\psi}^{p,q}_{C}(x)\left(\nabla \psi^{p}_{C}(x)-\nabla U(x)\right)p(x)dx\nonumber\\
&\qquad+\frac{1}{2\gamma}\left[\int_{E}\left|\nabla \psi_{C}^{p+\gamma h}(x)-\nabla U(x)\right|^{2}q(x)dx - \int_{E}\left|\nabla \psi_{C}^{p}(x)-\nabla U(x)\right|^{2}q(x)dx\right]\nonumber\\
&=\int_{E}\nabla\hat{\psi}^{p+\gamma h,q}_{C}(x)\left(\nabla \psi^{p+\gamma h}_{C}(x)-\nabla U(x)\right) h(x)dx+\int_{E}\nabla\hat{\psi}^{p+\gamma h,q}_{C}(x)\frac{\nabla \psi^{p+\gamma h}_{C}(x)-\nabla \psi^{p}_{C}(x)}{\gamma} p(x)dx\nonumber\\
&\qquad+\int_{E}\frac{\nabla \hat{\psi}^{p+\gamma h,q}_{C}(x)-\nabla \hat{\psi}^{p,q}_{C}(x)}{\gamma} \left(\nabla \psi^{p}_{C}(x)-\nabla U(x)\right)p(x)dx\nonumber\\
&\qquad+\frac{1}{2}\int_{E}\frac{\nabla \psi^{p+\gamma h}_{C}(x)-\nabla \psi^{p}_{C}(x)}{\gamma} \left(\nabla \psi_{C}^{p+\gamma h}(x)+\nabla \psi_{C}^{p}(x)-2\nabla U(x)\right)q(x)dx.
\end{align*}

As it was done for the computation of the first order directional derivative, we next notice that for every $g\in\mathcal{C}^{1}(E)$,  $\frac{\nabla \hat{\psi}^{p+\gamma h,q}_{C}(x)-\nabla \hat{\psi}^{p,q}_{C}(x)}{\gamma} $ satisfies
\begin{align*}
0&=\int_{E}\left[\frac{\nabla \hat{\psi}^{p+\gamma h,q}_{C}(x)-\nabla \hat{\psi}^{p,q}_{C}(x)}{\gamma}p(x)+\nabla \hat{\psi}^{p+\gamma h,q}_{C}(x)h(x)\right]\nabla g(x)dx\nonumber\\
&\qquad +\int_{E}\frac{\nabla \psi^{p+\gamma h}_{C}(x)-\nabla \psi^{p}_{C}(x)}{\gamma}q(x)\nabla g(x)dx.
\end{align*}

As in Section 3 of \cite{Gartner1977}, it follows then that there is a $\hat{\hat{\psi}}^{p,q,h}_{C}(x)\in \mathcal{C}^{(1+\alpha)}(E)$ such that
\begin{equation*}
\hat{\psi}^{p+\gamma h,q}_{C}(x)=\hat{\psi}^{p,q}_{C}(x)+\gamma \hat{\hat{\psi}}^{p,q,h}_{C}(x)+o_{2}(\gamma)
\end{equation*}
where $\left\Vert o_{2}(\gamma)\right\Vert_{(1+\alpha)}\rightarrow 0$ as $\gamma\rightarrow 0$. Then, for every $g\in\mathcal{C}^{1}(E)$, $\nabla\hat{\hat{\psi}}^{p,q,h}(x)$ satisfies
\begin{align}
\int_{E}\left[\nabla \hat{\hat{\psi}}^{p,q,h}_{C}(x)p(x)+\nabla \hat{\psi}^{p,q}_{C}(x)h(x)+\nabla \hat{\psi}^{p,h}_{C}(x)q(x)\right]\nabla g(x)dx&=0. \label{Eq:SecondDerivativeConstraint}
\end{align}

Let us then denote
\[
d^{2}J_{C}(p;q,h)=\lim_{\gamma\downarrow 0}\frac{dJ_{C}(p+\gamma h;q)-dJ_{C}(p;q)}{\gamma}.
\]

 We get
\begin{align}
d^{2}J_{C}(p;q,h)&=\int_{E}\left(\nabla\psi^{p,q}_{C}(x)-\nabla U(x)\right)\nabla \hat{\psi}^{p,q}_{C}(x)h(x)dx+\int_{E}\nabla\hat{\psi}^{p,q}_{C}(x)\nabla\hat{\psi}^{p,h}_{C}(x) p(x)dx\nonumber\\
&\qquad+\int_{E}\nabla \hat{\hat{\psi}}^{p, h,q}_{C}(x)\left(\nabla \psi^{p}_{C}(x)-\nabla U(x)\right)p(x)dx\nonumber\\
&\qquad+\int_{E}\nabla \hat{\psi}^{p,h}_{C}(x) \left(\nabla \psi_{C}^{p}(x)-\nabla U(x)\right)q(x)dx \nonumber\\
&=\int_{E}\left(\nabla \psi_{C}^{p}(x)-\nabla U(x)\right)\left( \nabla\hat{\psi}^{p,q}(x)h(x)+\nabla\hat{\psi}^{p,h}(x)q(x)+\nabla\hat{\hat{\psi}}^{p,q,h}(x)p(x)\right)dx\nonumber\\
&\qquad +\int_{E}\nabla\hat{\psi}^{p,q}(x)\nabla\hat{\psi}^{p,h}(x)p(x)dx.\label{Eq:SecondDerivativeRateFcn}
\end{align}

Using the constraint (\ref{Eq:SecondDerivativeConstraint}) with the test function $g(x)=\psi_{C}^{p}(x)-U(x)$, we then obtain
\begin{align}
d^{2}J_{C}(p;q,h)
&=\int_{E}\nabla\hat{\psi}^{p,q}(x)\nabla\hat{\psi}^{p,h}(x)p(x)dx.\nonumber 
\end{align}

Recall that for every  $g\in\mathcal{C}^{1}(E)$, $\nabla\hat{\psi}^{p,q}(x)$ satisfies  (\ref{Eq:FirstDerivativeConstraint}) and similarly for $\nabla\hat{\psi}^{p,h}(x)$. Thus, selecting $h(x)=q(x)$, we get
\begin{align}
d^{2}J_{C}(p;q,q)
&=\int_{E}\left|\nabla\hat{\psi}^{p,q}_{C}(x)\right|^{2}p(x) dx\label{Eq:SecondDerivativeRateFcnInvMeas2}
\end{align}

 Relation (\ref{Eq:SecondDerivativeRateFcnInvMeas2}) implies that pointwise in $p$ and for non-zero directions $q(x)$ the second order directional derivative of $I_{C}(p)$ increases when adding an appropriately non-zero irreversible drift $C$, i.e.,
\[
d^{2}J_{C}(p;q,q)\geq 0.
\]

Of course, this is expected to be true due to convexity. Let us next investigate what happens at the law of large numbers limit $\mu=\bar{\pi}$. So, let us choose $\mu(dx)$ to be the invariant measure $\bar{\pi}(dx)$ and let us denote its density by $\bar{p}(x)$. Then, we notice that in this case $\nabla\psi^{\bar{p}}_{C}(x)=\nabla U(x)$. So, (\ref{Eq:SecondDerivativeRateFcnInvMeas2}) becomes
\begin{align}
d^{2}J_{C}(\bar{p};q,q)
&=\int_{E}\left|\nabla\hat{\psi}^{\bar{p},q}_{C}(x)\right|^{2}\bar{\pi}(dx)\geq 0 \label{Eq:SecondDerivativeRateFcnInvMeas3}
\end{align}
where, $\forall g\in\mathcal{C}^{1}(E)$, $\nabla\hat{\psi}^{\bar{p},q}(x)$ satisfies
\begin{align}
\int_{E}\left[\nabla \hat{\psi}^{\bar{p},q}_{C}(x) \bar{p}(x)+C(x)q(x)\right]\nabla g(x)dx=0. \label{Eq:FirstDerivativeConstraintInvMeas2}
\end{align}

In fact, we get for $q,C$ such that $\textrm{div}(qC)\neq 0$  that $\nabla \hat{\psi}^{p,q}_{C}(x)\neq 0$. Then, by (\ref{Eq:FirstDerivativeConstraintInvMeas2}) and (\ref{Eq:SecondDerivativeRateFcnInvMeas3}) we have
\begin{equation}
d^{2}J_{C}(\bar{p};q,q)>0\label{Eq:SecondGateauxDerivative}.
\end{equation}

In addition, (\ref{Eq:FirstDerivativeConstraintInvMeas2}) shows that if $C=0$, or if $q,C$ are such that $\text{div}(qC)=0$, then
\begin{equation*}
d^{2}J_{0}(\bar{p};q,q)=0.
\end{equation*}

\subsection{Completion of the proof of Theorem \ref{T:VarianceObservable}}\label{SSS:ProofAsymptoticVariance}

Let $C \not=0$ and $f \in \mathcal{C}^{(\alpha)}(E)$ be such that $\tilde{I}_{f,C}(\ell)>\tilde{I}_{f,0}(\ell)$ for every $\ell\neq\bar{f}$. Then, we want to prove
\begin{equation*}
\tilde{I}_{f,C}^{\prime\prime}(\bar{f})-\tilde{I}_{f,0}^{\prime\prime}(\bar{f})>0.
\end{equation*}

We know by Proposition \ref{P:ExistenceMinimizer} that there exist measures, say $\mu_{C}(dx;\ell)$ and  $\mu_{0}(dx;\ell)$, that have a strictly positive densities in $\mathcal{C}^{(2+\alpha)}(E)$ such that
\begin{equation*}
\tilde{I}_{f,C}(\ell)=I_{C}(\mu_{C}(\cdot; \ell)),\quad\text{ and }\quad  \tilde{I}_{f,0}(\ell)=I_{0}(\mu_{0}(\cdot; \ell))
\end{equation*}

By convexity and the definitions of $\mu_{C}(dx;\ell)=p_{C,\ell}(x)dx$ and  $\mu_{0}(dx;\ell)=p_{0,\ell}(x)dx$, we have that for all $\ell\in\left(\min_{x}f(x),\max_{x}f(x)\right)$
\begin{align*}
\tilde{I}_{f,C}^{\prime\prime}(\ell)-\tilde{I}_{f,0}^{\prime\prime}(\ell)&=\frac{\partial^{2}}{\partial \ell^{2}}\left[I_{C}(\mu_{C}(\cdot; \ell))-I_{0}(\mu_{0}(\cdot; \ell))\right]\nonumber\\
&=\frac{\partial^{2}}{\partial \ell^{2}}\left[ \left(I_{C}\left(\mu_{C}(\cdot; \ell)\right)-I_{0}\left(\mu_{C}(\cdot; \ell)\right)\right)+\left(I_{0}(\mu_{C}(\cdot; \ell))-I_{0}(\mu_{0}(\cdot; \ell))\right)\right]\nonumber\\
&\geq \frac{\partial^{2}}{\partial \ell^{2}}\left[I_{C}(\mu_{C}(\cdot; \ell))-I_{0}(\mu_{C}(\cdot; \ell))\right]\nonumber\\
&= \frac{\partial^{2}}{\partial \ell^{2}} J_{C}(\mu_{C}(\cdot; \ell)).\nonumber\\
\end{align*}

Then, (\ref{Eq:SecondDerivativeRateFcnInvMeas3}) implies that when evaluated at the law of large numbers $\ell=\bar{f}$,
\begin{align}
\frac{\partial^{2}}{\partial \ell^{2}} J_{C}(p_{C,\ell})\Big|_{\ell=\bar{f}}
&=\int_{E}\left|\nabla\hat{\psi}^{\bar{p},\bar{q}}_{C}(x)\right|^{2}\bar{\pi}(dx),\label{Eq:VarProblem1}
\end{align}
such that (\ref{Eq:FirstDerivativeConstraintInvMeas2}) holds with $q(x)=\bar{q}(x)=\frac{\partial}{\partial\ell}p_{C}(x;\ell)\Big |_{\ell=\bar{f}}$, i.e.
 $\nabla\hat{\psi}^{\bar{p},q}(x)$ satisfies
\begin{align}
\int_{E}\left[\nabla \hat{\psi}^{\bar{p},\bar{q}}_{C}(x) \bar{p}(x)+C(x)\bar{q}(x)\right]\nabla g(x)dx=0,\quad \forall g\in\mathcal{C}^{1}(E). \label{Eq:VarProblem2}
\end{align}

Then, (\ref{Eq:SecondGateauxDerivative}) implies
\begin{align*}
\frac{\partial^{2}}{\partial \ell^{2}} J_{C}(p_{C,\ell})\Big|_{\ell=\bar{f}}
&>0,
\end{align*}
as long as $\textrm{div}\left(\bar{q}C\right)\neq 0$. This concludes the proof of Theorem \ref{T:VarianceObservable}.

\qed
\section{Simulations}\label{S:Simulations}

In this section  we present some numerical results to illustrate the theoretical findings. We study numerically the effect that adding irreversibility has on the speed of convergence to the equilibrium. Consider the SDE in 2 dimensions
\begin{eqnarray*}
dZ_{t}=\left[-\nabla U(Z_{t})+C(Z_{t})\right]dt+\sqrt{2D}dW_{t}, \quad Z_{0}=0
\end{eqnarray*}
where $D=0.1$ and, for $z=(x,y)$, $C(x,y)=\delta C_{0}(x,y)$ with $C_{0}(x,y)=J \nabla U(x,y)$. Here, $\delta\in \mathbb{R}$, $I$ is the $2\times2$ identity matrix and $J$ is the standard $2\times2$ antisymmetric matrix, i.e., $J_{12}=1$, $J_{21}=-1$ and $J_{11}=J_{22}=0$.

Clearly, in the case $\delta=0$ we have reversible dynamics, whereas for $\delta\neq 0$ the dynamics is irreversible. Notice that for any $\delta\in\mathbb{R}$, the invariant measure is
\begin{equation*}
\bar{\pi}(dxdy)=\frac{e^{-\frac{U(x,y)}{D}}}{\int_{\mathbb{R}^{2}}e^{-\frac{U(x,y)}{D}}dxdy}dxdy
\end{equation*}

Let us suppose that we are given an observable $f(x,y)$ and we want to compute
\begin{equation*}
\bar{f}=\int_{\mathbb{R}^{2}}f\left(x,y\right)\bar{\pi}(dxdy).
\end{equation*}

It is known that an estimator for $\bar{f}$ is given by
\begin{equation*}
\hat{\bar{f}}(t)=\frac{1}{t-v}\int_{v}^{t}f\left(X_{s},Y_{s}\right)ds
\end{equation*}
where $v$ is some burn-in period that is used with the hope that the bias has been significantly reduced by time $v$. This estimate is based on simulating a very long trajectory $Z_{s}=\left(X_{s},Y_{s}\right)$.

In general, a central limit theorem holds and takes the following form
\[
t^{1/2}\left(\hat{\bar{f}}(t)-\bar{f}\right)\Rightarrow N(0,\sigma^{2}_{f})
\]
where $\sigma^{2}_{f}$ is the asymptotic variance and is a deterministic constant. Then it is known, e.g., Proposition IV.1.3 in \cite{AsmussenGlynn2007}, that
\begin{equation*}
\sigma^{2}_{f}=2\int_{0}^{\infty}c(s)ds, \quad \text{where}\quad c(s)=\mathbb{E}_{\bar{\pi}}\left[\left(f(X_{0},Y_{0})-\bar{f}\right)\left(f(X_{s},Y_{s})-\bar{f}\right)\right]
\end{equation*}

The objective now is to see how $\sigma^{2}_{f}$ scales as a function of $\delta$. For this purpose, we recall that  up to constants
\begin{equation*}
\sigma^{2}_{f}=\frac{1}{2 \tilde{I}^{\prime\prime}(\bar{f})}
\end{equation*}
where $\tilde{I}^{\prime\prime}(\bar{f})$ is the second derivative of the large deviations action functional evaluated at $\ell=\bar{f}$.

We have seen already that adding the irreversibility $C$ in the dynamics results in  smaller variance $\sigma^{2}_{f}$, as Theorem \ref{T:VarianceObservable} verifies. Let us demonstrate this through an empirical study. To do so we use the well established method of batch means (e.g., Section IV.5 in \cite{AsmussenGlynn2007}) in order to construct confidence interval for $\hat{\bar{f}}(t)$. Let us recall here the algorithm for convenience.

Let us fix a desired time instance $t$ and the number of batches, say $m$. Then for $\kappa=1,\cdots,m$ we define
\[
\hat{\bar{f}}(t;\kappa)=\frac{1}{t/m}\int_{(\kappa-1)t/m}^{\kappa t/m}f\left(X_{s},Y_{s}\right)ds,
\]
\[
\hat{\bar{f}}(t)=\frac{1}{m}\sum_{\kappa=1}^{m}\hat{\bar{f}}(t;\kappa)
\]
and
\[
s^{2}_{m}(t)=\frac{1}{m-1}\sum_{\kappa=1}^{m}\left(\hat{\bar{f}}(t;\kappa)-\hat{\bar{f}}(t)\right)^{2}
\]

Then, we have in distribution
\[
\sqrt{m}\frac{\hat{\bar{f}}(t)-\bar{f}}{s_{m}(t)}\Rightarrow T_{m-1}, \quad\text{ as }t\rightarrow\infty
\]
where $T_{m-1}$ is the Student's T distribution with $m-1$ degrees of freedom. So, a $(1-\alpha)\%$ confidence interval is given by
\[
\left(\hat{\bar{f}}(t)-t_{\alpha/2,m-1}s_{m}(t)/\sqrt{m},\hat{\bar{f}}(t)+t_{\alpha/2,m-1}s_{m}(t)/\sqrt{m}\right)
\]

For the simulations that follow we used time step $\Delta t=0.001$,\ and number of batches ranging from $m=10$ to $m=20$ at $t$ gets larger. Also, in order to minimize the bias, we used  a burn-in time $v=5$.

\begin{figure}
[ptb]
\begin{center}
\includegraphics[height=6cm, width=12cm]
{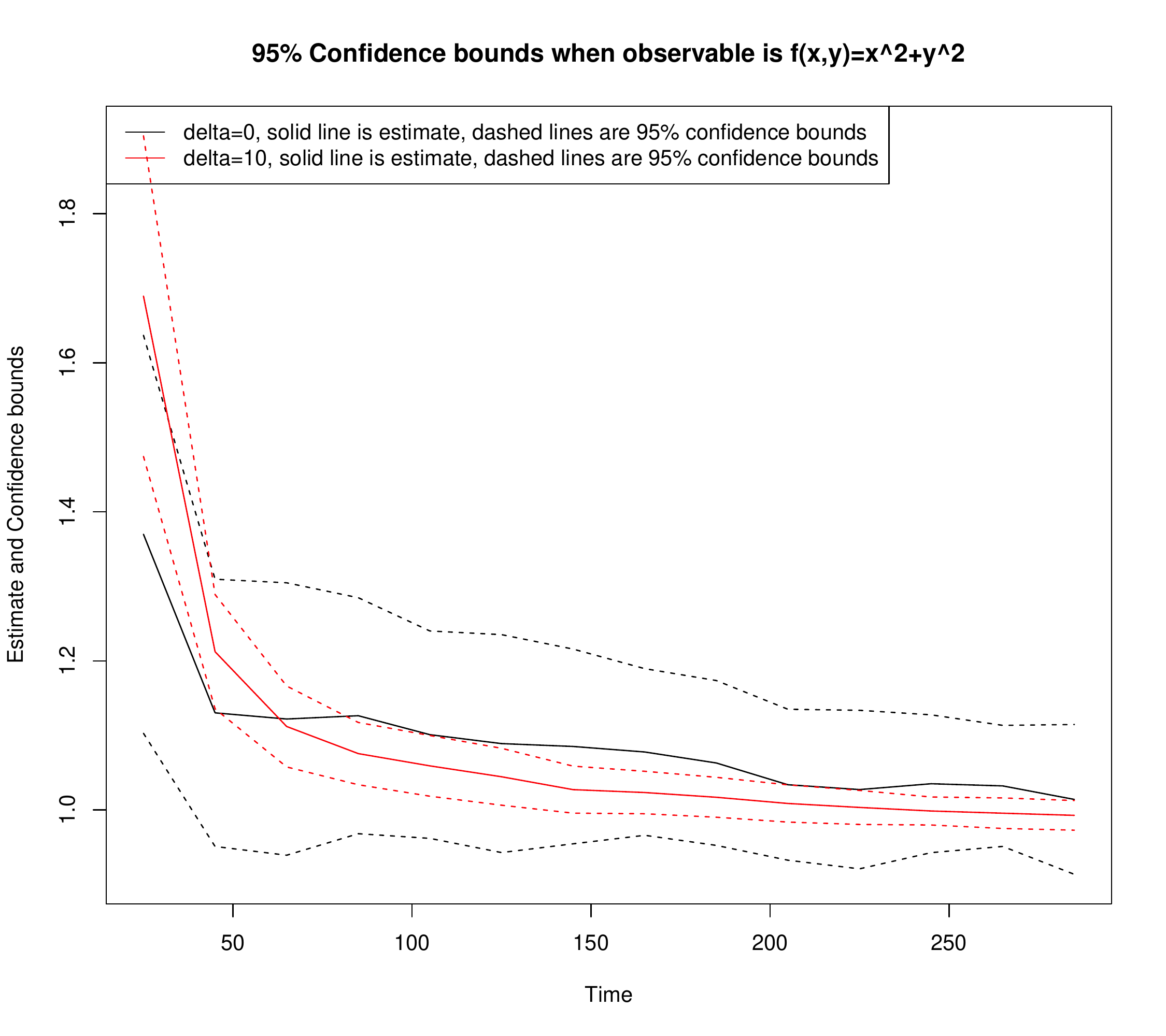}%
\caption{Estimate and $95\%$ Confidence bounds when $U(x,y)=\frac{1}{4}(x^{2}-1)^{2}+\frac{1}{2}y^{2}$ and $f(x,y)=x^2 + y^2$.}%
\label{Fig3}%
\end{center}
\end{figure}

\begin{figure}
[ptb]
\begin{center}
\includegraphics[height=6cm, width=12cm]
{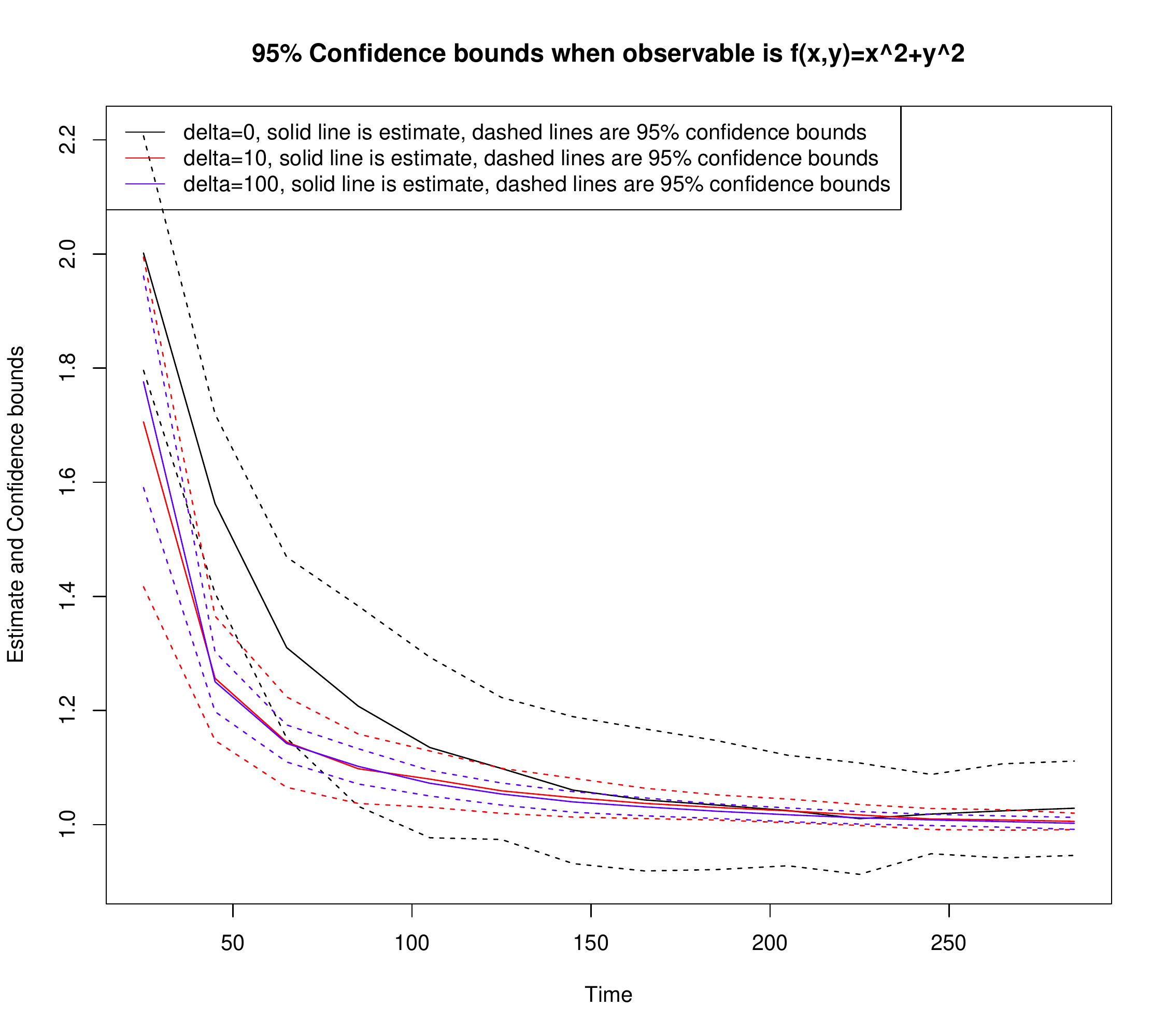}%
\caption{Estimate and $95\%$ Confidence bounds when $U(x,y)=\frac{1}{4}(x^{2}-1)^{2}+\frac{1}{2}y^{2}$ and $f(x,y)=x^2 + y^2$.}%
\label{Fig4}%
\end{center}
\end{figure}

We present three different examples. In the first example we pick the potential $U(x,y)=\frac{1}{4}(x^{2}-1)^{2}+\frac{1}{2}y^{2}$ and the observable $f(x,y)=x^{2}+y^{2}$. These dynamics was also considered in \cite{LelievreNierPavliotis2012}. We remark here that the quantity $\int_{E}\left(x^{2}+y^{2}\right)\bar{\pi}(dx,dy)$ is the long-time mean-square displacement of the process $Z_{t}=(X_{t},Y_{t})$. In Figures \ref{Fig3} and \ref{Fig4} we see $95\%$ confidence bounds for $\hat{\bar{f}}(t)$. It is clear the adding irreversibility not only speeds up convergence to equilibrium, but it also results in significant reduction in the variance. In Figure \ref{Fig3} we compare the reversible case (i.e., with $\delta=0$) with the irreversible case with $\delta=10$. Then, in Figure \ref{Fig4}, we have also included the case $\delta=100$. For the particular test case, the confidence bounds are even tighter when $\delta=100$ when compared to $\delta=10$. This result illustrates Theorem \ref{T:measure2} and Theorem \ref{T:VarianceObservable}.

For illustration purposes, we present in Table \ref{Table1a1}, variance estimates for different values of $\delta $ and time horizons $t$ in the set-up of Figure \ref{Fig4}. It is noteworthy that the variance reduction for this particular example is about two orders of magnitude.

\begin{table}[th]
\begin{center}
{\small
\begin{tabular}
[c]{|c|c|c|c|c|c|}\hline
$\delta\hspace{0.1cm} | \hspace{0.1cm} t$ & $25$ & $100$ & $160$ & $220$ &
$295$ \\\hline
$0$ & $0.22$ & $0.08$ & $0.038$ & $0.029$ & $0.011$ \\\hline
$10$ & $0.19$ & $0.01$ & $0.007$ & $0.005$ & $0.002$ \\\hline
$100$ & $0.09$ & $0.001$ & $3e-04$ & $2.8e-04$ & $1.3e-04$ \\\hline
\end{tabular}
}
\end{center}
\caption{Estimated variance values for different pairs $(\delta,t)$.}%
\label{Table1a1}%
\end{table}

\begin{figure}[ptb]
\includegraphics[height=6cm,width=12cm]{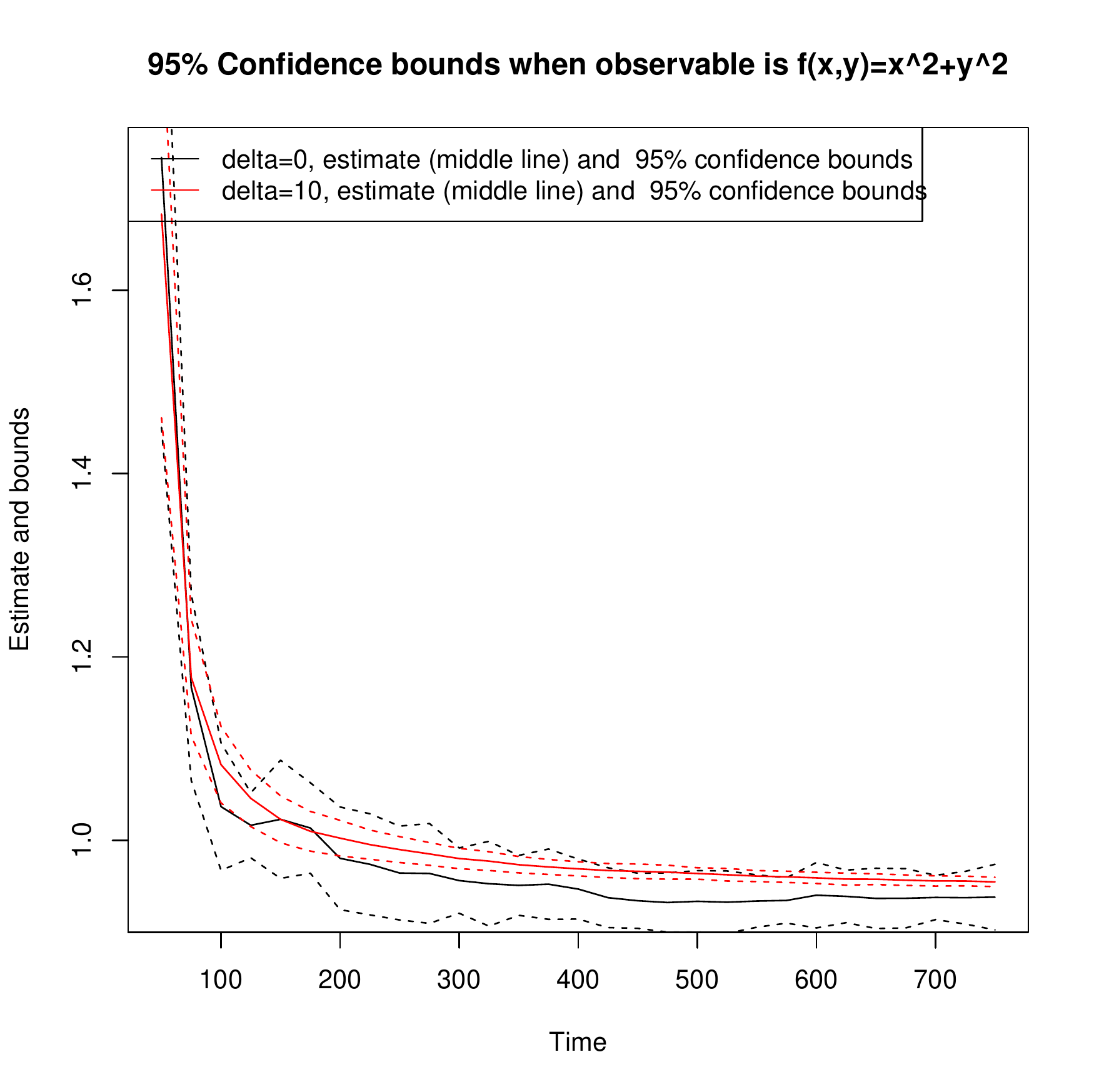}%
\caption{Estimate and $95\%$ Confidence bounds when $U(x,y)=(x^{2}-1)^{2}+\frac{1}{2}(3y+x^{2}-1)^{2}$ and $f(x,y)=x^2 + y^2$.} %
\label{Fig2}%
\end{figure}

In the second example we pick again a bimodal potential $U(x,y)=(x^{2}-1)^{2}+\frac{1}{2}(3y+x^{2}-1)^{2}$ and the observable $f(x,y)=x^{2}+y^{2}$. In Figure \ref{Fig2} we see $95\%$ confidence bounds for $\hat{\bar{f}}(t)$. In Table \ref{Table2a1}, we present numerical data for the variance estimates that are illustrated in Figure \ref{Fig2}. Again, we see variance reduction and it is at the order of about two  magnitudes.

\begin{table}[th]
\begin{center}
{\small
\begin{tabular}
[c]{|c|c|c|c|c|c|c|c|}\hline
$\delta\hspace{0.1cm} | \hspace{0.1cm} t$ & $100$ & $200$ & $300$ & $400$ & $500$ &
$600$ & $700$ \\\hline
$0$ &  $0.01$ & $0.006$ & $0.002$ & $0.002$ & $0.002$ & $0.003 $ & $0.002 $\\\hline
$10$ &  $0.003$ & $0.0007$ & $0.0002$ & $0.0001$ & $7e-05$ & $6e-05 $ & $6e-05 $\\\hline
\end{tabular}
}
\end{center}
\caption{Estimated variance values for different pairs $(\delta,t)$.}%
\label{Table2a1}%
\end{table}

In the third example we pick the potential
\[
U(x,y)=\frac{1}{4}
\Big[(x^2-1)^2((y^2-2)^2+1)+2y^2 - y/8\Big]+ e^{-8x^{2}-4y^{2}}.
\]

\begin{figure}[ptb]
\includegraphics[height=4cm,width=8cm]{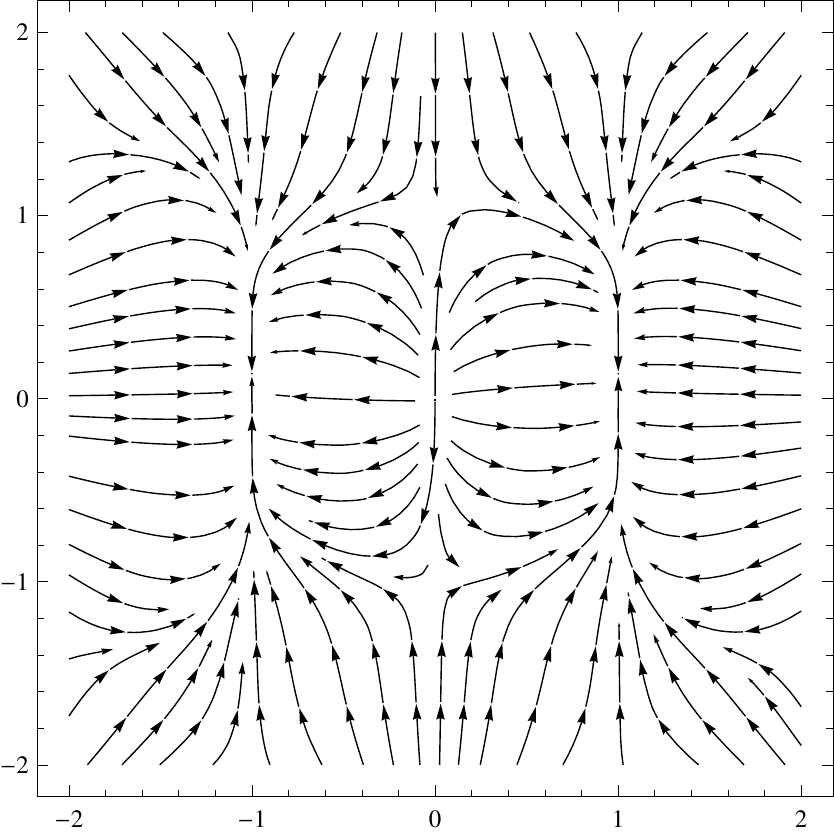}%
\caption{Phase portrait of $U(x,y)=\frac{1}{4}
\Big[(x^2-1)^2((y^2-2)^2+1)+2y^2 - y/8\Big]+ e^{-8x^{2}-4y^{2}}$.} %
\label{F:PhasePortrait}%
\end{figure}

Due to the somewhat complex form of $U(x,y)$, we have also plotted  in Figure \ref{F:PhasePortrait} its phase portrait.
We see that it has two local minima at $(\pm 1.00051, 0.125314)$, two saddle points at $(0,-1.00711)$ and at $(0, 1.08849)$ and a local maximum at $(0,-0.0139)$.

\begin{figure}[ptb]
\includegraphics[height=12cm,width=6cm, angle=-90]{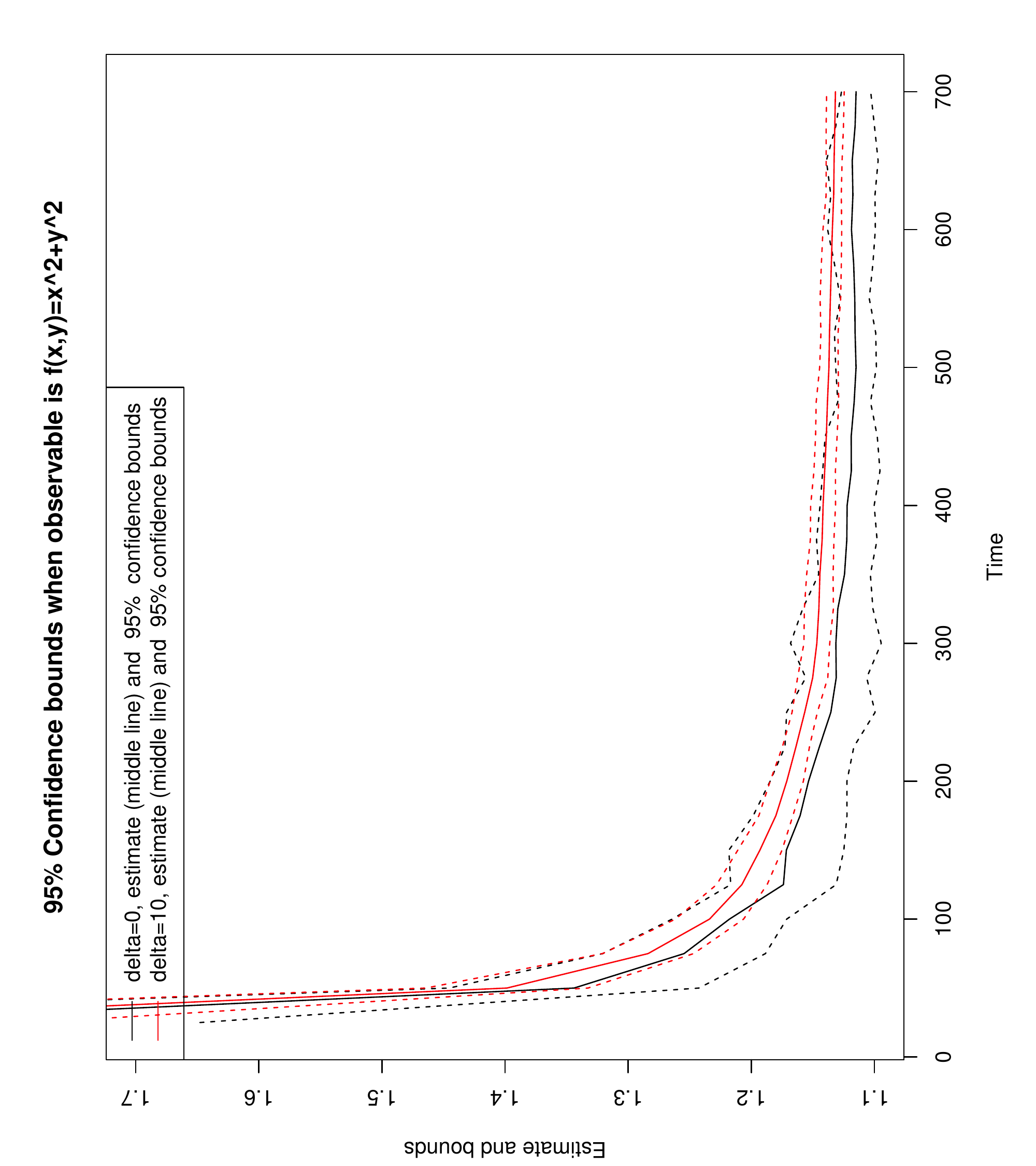}%
\caption{Estimate and $95\%$ Confidence bounds when $U(x,y)=\frac{1}{4}
\Big[(x^2-1)^2((y^2-2)^2+1)+2y^2 - y/8\Big]+ e^{-8x^{2}-4y^{2}}$ and $f(x,y)=x^2 + y^2$.} %
\label{Fig5}%
\end{figure}

We consider again the observable $f(x,y)=x^{2}+y^{2}$. In Figure \ref{Fig5} we see $95\%$ confidence bounds for $\hat{\bar{f}}(t)$. In Table \ref{Table3a1}, we present numerical data for the variance estimates that are illustrated in Figure \ref{Fig5}. Again, we see variance reduction and it is at the order of about one  magnitude when the irreversible parameter is $\delta=10$.

\begin{table}[th]
\begin{center}
{\small
\begin{tabular}
[c]{|c|c|c|c|c|c|c|c|c|}\hline
$\delta\hspace{0.1cm} | \hspace{0.1cm} t$ & $100$ & $200$ & $300$ & $400$ & $500$ & $600$  & $700$  \\\hline
$0$  &  $0.004$ & $0.002$ & $0.002$ & $0.001$ & $0.001$ & $0.001 $ & $0.001$\\\hline
$10$ &  $0.001$ & $0.0003$ & $0.0002$ & $0.0001$ & $0.0001$ & $0.0001 $ &$0.0001$\\\hline
\end{tabular}
}
\end{center}
\caption{Estimated variance values for different pairs $(\delta,t)$.}%
\label{Table3a1}%
\end{table}

We conclude this section with a remark on the optimal choice of irreversibility. Theorem \ref{T:measure2} suggests that in the generic situation, perturbations of the form $C(\cdot)=\delta C_{0}(\cdot)$ yield better results as the parameter $\delta$ increases. However, in practice the higher the $\delta$ is, the smaller the discretization step in the simulation algorithm should be, i.e., there is a trade-off to consider here. Thus it makes sense to look for the optimal perturbation $C(x)$ and this could be formulated as a solution to a variational problem that involves minimizing the asymptotic variance of the estimator. Since, the asymptotic variance is inversely proportional to the second derivative of the rate function of the observable evaluated at $\bar{f}$, the variational problem to consider is basically maximization over vector fields $C$ that satisfy condition $\bf{(H)}$ of the quantity (\ref{Eq:VarProblem1}) under the constraint  (\ref{Eq:VarProblem2}). We plan to investigate this question in a
future work.

\section{Conclusions}\label{S:Conclusions}

In this article we have considered the problem of estimating the expected value of a functional of interest using as estimator the long time average of a process that has as its invariant distribution the target measure. We have argued using large deviations theory, both theoretically and numerically,  that adding an appropriate drift to the dynamics of a reversible Langevin equation, results in smaller asymptotic variance for the time average estimator.  We  characterize when observables do not see their variance reduced in terms of a precise non-linear Poisson equation.

\end{document}